\RequirePackage{tikz}
\documentclass[a4paper,11pt]{amsart}

\usepackage{graphicx}
\usepackage{amssymb}
\usepackage{amsmath}
\usepackage{amsfonts}
\usepackage{bm}
\usepackage{xcolor}
\usepackage{tabularx}
\usepackage{subcaption}

\graphicspath{{images/}, {images/non_parametric/}, {images/geometrical_parametrization/}, {images/physical_geometrical_parametrization/}}

\usepackage[foot]{amsaddr}
\usepackage[margin=2.7cm]{geometry}
\theoremstyle{thmstyletwo}%
\newtheorem{remark}{Remark}%
\newtheorem{definition}{Definition}%

\begin{document}

\title[Partitioned ROMs for parametrized FSI problems]{Projection based semi--implicit partitioned Reduced Basis Method for Fluid--Structure Interaction problems}
\author{Monica Nonino\textsuperscript{1}}
\email{monica.nonino@univie.ac.at}

\author{Francesco Ballarin\textsuperscript{2}}
\email{fballarin@unicatt.it}

\author{Gianluigi Rozza\textsuperscript{3}}
\email{grozza@sissa.it}

\author{Yvon Maday\textsuperscript{4}}
\email{maday@ann.jussieu.fr}

\address{\textsuperscript{1}University of Vienna, Department of Mathematics.}
\address{\textsuperscript{2}Universit\'a Cattolica del Sacro Cuore, Dipartimento di Matematica e Fisica.}
\address{\textsuperscript{3}Sissa, mathematics area, mathLab, International School for Advanced Studies.}
\address{\textsuperscript{4}Laboratoire Jacques-Louis Lions, Universit\'e Pierre et Marie Curie.}

\subjclass[2010]{78M34, 97N40, 35Q35}
\begin{abstract}
In this manuscript a POD--Galerkin based Reduced Order Model for unsteady Fluid--Structure Interaction problems is presented. The model is based on a partitioned algorithm, with semi--implicit treatment of the coupling conditions. A Chorin--Temam projection scheme is applied to the incompressible Navier--Stokes problem, and a Robin coupling condition is used for the coupling between the fluid and the solid. The coupled problem is based on an Arbitrary Lagrangian Eulerian formulation, and the Proper Orthogonal Decomposition procedure is used for the generation of the reduced basis. We extend existing works on a segregated Reduced Order Model for Fluid--Structure Interaction to unsteady problems that couple an incompressible, Newtonian fluid with a linear elastic solid, in two spatial dimensions. We consider three test cases to assess the overall capabilities of the method: an unsteady, non--parametrized problem, a problem that presents a geometrical parametrization of the solid domain, and finally, a problem where a parametrization of the solid's shear modulus is taken into account.
\end{abstract}

\maketitle

\section{Introduction}
\noindent Fluid--Structure Interaction (FSI) problems are a wide spread topic in the applied mathematics community, and despite their instrinsic complicated nature (see for example \cite{MadayCRAS, GRANDMONT1998525}), they are frequently used for simulation purposes, for example, in naval engineering \cite{LombardiParolini}, as well as in biomedical applications (as an example of the implementation of FSI in the medical field see \cite{QuainiWang2017, Quaini:129770, BallarinJCP, String_model_1, String_model_2, QuainiWang2017, Maday2009, BertagnaVeneziani}) and in aeronautical engineering (see for example \cite{FarhatPeterson, FarhatAvery, FARHAT201453, LIEU20065730, AmsallemCortialFarhat, Carlberg2013623}).\\
The complex nature of these problems is reflected not only by their theoretical treatment, but also by the way they are treated numerically.
There are two approaches that can be adopted in order to address a FSI problem: the first approach consists of a monolithic procedure \cite{BallarinRozza, NoBaRo19, BADIA20097986}, whereas the second approach consist of a partitioned (segregated), procedure \cite{Ballarin2017, Fernandez, kalashnikova2013stable}. 
\\In a monolithic algorithm the fluid and the solid problem are solved \emph{simultaneously}: this results in algorithms that {show good stability properties, with respect to time, independently of the density of the fluid and the solid, and independently of the geometrical properties of the physical domain; indeed the monolithic algorithm does not suffer from the so called \emph{added mass effect} (see \cite{CausinGerbeauNobile} for an analytical study of this phenomenon), which is very well known in the FSI community, and is responsible for the numerical instabilities in the design of partitioned algorithms. Stability in time} is highly desirable in the framework of unsteady problems, especially if we wish to use large time--steps in the simulations; the main drawback of these monolithic algorithms is given by the fact that they deeply rely on the availability of legacy softwares that can be used to solve both the fluid problem and the solid problem: in this sense, monolithic algorithms are less flexible and more tailored to the particular problem at hand. In the literature there are many examples of works that are based on a monolithic approach: in \cite{Quaini2017} the author focuses on a monolithic algorithm to address a coupled problem, written within the Arbitrary Lagrangian Eulerian (ALE) formalism, which models the interaction between the blood flow and the arterial walls; another example of FSI problems related to the blood flow--arterial interaction can be found in \cite{WuCai}. In \cite{GeeWall} the authors propose different preconditioners, to be used in a Newton--Krylov method for the nonlinear problem arising from solving in a monolithic fashion a coupled problem. For the reader interested in a general introduction to monolithic approaches to FSI problems, we refer to \cite{Richter}.
\\As an alternative to monolithic approaches, one can think of adopting a partitioned procedure; indeed, existing simulation tools for fluid dynamics and for structural dynamics are well developed and are used on a daily basis in industrial applications. It is therefore natural to try to combine these computational tools, to address coupled problems: this is exactly the rationale behind a partitioned algorithm. In a partitioned procedure, we solve \emph{separately} the fluid and the solid problem, and then we couple the two physics with some iterative procedure, see for example \cite{GuimetMaday}. 
Partitioned approaches are very flexible, as they allow to design the procedure in different ways, according to the problem under consideration. In \cite{Cesmelioglu2016}, the authors propose a segregated algorithm to solve a FSI problem, where the coupling of the two physics at the fluid--structure interface is taken care of through a constrained optimization problem. In \cite{FernandezVidrascu, FernandezMullaertVidrascu} the authors consider the problem of coupling an incompressible fluid with a thin structure; in \cite{FernandezVidrascu} the authors propose a Robin--Neumann type of coupling at the fluid--structure interface, whereas in \cite{FernandezMullaertVidrascu} the authors propose and explain different couplings techniques at the  fluid--structure interface, within an explicit coupling setting. On the contrary, in \cite{FernandezGerbeauGrandmont} the authors deal with a problem that has a strong added mass effect, which is typically the case for the blood in the vessels: here, an implicit coupling is the way to control the stability issues due to the added mass effect. Nevertheless, it is clear that a fully implicit treatment of the coupling conditions leads to prohibitive computational costs; for this reason, in \cite{FernandezGerbeauGrandmont} the authors propose a semi--implicit coupling technique, which is the approach that will be adopted in this manuscript.
\\Addressing a coupled problem by means of a partitioned procedure is advantageous in terms of computational efficiency, also from the Reduced Order Model (ROM) point of view: indeed, in the online phase of the Reduced Basis Method \cite{RoHuMa13, Rozza2008229, RoVe07, LassilaQuarteroniRozza, Lassila, osti_974411, HesthavenRozzaStamm}, we have to solve, separately, smaller systems. Moreover, with some minor changes such as change of variables and appropriate choices for the couplings, it is possible to further reduce the dimension of the online systems, as we will see in the following. {In the model order reduction framework, there is also a fair amount of work that is being carried out and which focuses on ROM--ROM and ROM--FOM coupling, see for example \cite{HOANG2021113997, TezaurKuberry, BERGMANN2018301}. All these works represent an extremely interesting approach from which many FSI applications of interest could benefit; for this reason, the authors believe that this direction represents a future line of work within partitioned algorithms.}
\\In this manuscript we design a segregated procedure, combined with a Reduced Order Model based on a Proper Orthogonal Decomposition. The goal is to extend the work done in \cite{Ballarin2017}, moving to the treatment of a two dimensional structure within an Arbitrary Lagrangian Eulerian formalism, and the work done in \cite{AstorinoChoulyFernandez, BadiaNobileVergara}, adapting the computation and the treatment of the Robin coupling condition, also to the case of a thick, two dimensional structure. {The present manuscript represents also an extension of the work done in \cite{monica}, where the problem under consideration was only unsteady, but no geometric or physical parametrization has been considered}.
\\This manuscript is structured as follows: in Section \ref{ALE formalism} we briefly introduce the Arbitrary Lagrangian Formulation, and we set the notation that will be used throughout the manuscript.
In Section \ref{problem formulation} we introduce the first test case, namely a time dependent, non parametrized FSI problem that models the interaction of a fluid with a thick, two dimensional, structure; in Section \ref{offline phase} we introduce the partitioned procedure at the high order level. In Section \ref{online phase} we derive the partitioned procedure at the reduced order level, and in Section \ref{numerical section non parametrized} we present the numerical results. 
In Section \ref{geometrical parametrization} we consider the same problem of interest, with the addition of a shape parametrization: in Section \ref{ALE geometrical formalism} we present the ALE formalism in the presence of a geometrical parametrization of the domain; in Section \ref{parametric problem formulation} we give the strong formulation of the problem of interest, and in Section \ref{offline phase parametrized} we describe the algorithm at the high order level. In Section \ref{shape and physical parametrization: online phase} we introduce the reduced order model, and then we present some numerical results in Section \ref{numerical results parametrized}, for the geometrical parametrization only. Then, in Section \ref{shape and physical parametrization: numerical results} we show some numerical results also in the presence of a physical parameter. Conclusions and considerations on future possible lines of work are presented in Section \ref{conclusions}.
\section{Configurations, definitions and notation}\label{ALE formalism}
In this section we are going to introduce briefly the Arbitrary Lagrangian Eulerian (ALE) formalism, in order to set the notation that will be used throughout the rest of this manuscript.\\
In FSI problems the fluid domain is a moving domain (except for those situations in which the displacement of the solid is very small, and thus the whole physical domain can be considered as fixed). In solid mechanics, on the other hand, it is common to deal with deforming domains, and the deformation itself is the unknown of the problem; for fluid dynamics instead one usually considers fixed domains. This different point of view is the motivation behind a formalism, very known and widely used in the community, which is called the \emph{Arbitray Lagrangian Eulerian formulation} \cite{Richter, Donea, Hughes, QuainiBasting2017}.
\\Let $\Omega(t)\subset\mathbb{R}^2$ be the physical domain over which the FSI problem is formulated, with time $t\in[0, T]$: $\Omega(t) = \Omega_f(t)\cup\Omega_s(t)$, where $\Omega_f(t)\subset\mathbb{R}^2$ and $\Omega_s(t)\subset\mathbb{R}^2$ are the fluid and the solid domain at time $t$, respectively; we assume that the two domains do not overlap, i.e. $\Omega_f(t)\cap\Omega_s(t)=\emptyset$, and finally, the fluid--structure interface $\Gamma_{FSI}(t)$ is defined as $\Gamma_{FSI}(t):=\bar{\Omega}_f(t)\cap\bar{\Omega}_s(t)$. 
To describe the behavior of a solid it is common practice to use the so called \emph{Lagrangian formalism}: all the quantities and the conservation laws are formulated on the reference configuration $\hat{\Omega}_s=\Omega_s(t=0)$. On the contrary, when describing the behavior of a fluid, the \emph{Eulerian formalism} is used instead: all the quantities and the conservation laws are formulated on the configuration $\Omega_f(t)$ at the current time $t$. In order to be able to describe both the fluid and the solid, a mixed formulation (the ALE formulation indeed) is used: the underlying idea is that of pulling back the fluid equations to an arbitrary time--independent configuration $\hat{\Omega}_f$: one possible choice for $\hat{\Omega}_f$ is $\hat{\Omega}_f = \Omega_f(t=0)$, the domain at initial time.
\begin{figure}
\centering
\begin{tikzpicture}
\node[anchor=south west,inner sep=0] (image) at (0,0) {\includegraphics[scale=0.2]{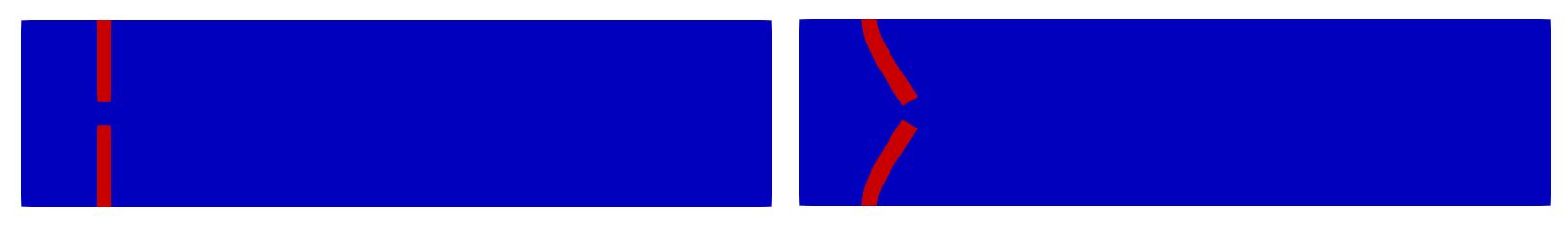}};
\begin{scope}[x={(image.south east)},y={(image.north west)}]
	\draw[white,thick] (0.9,0.5) node {$\Omega(t)$};
	\draw[white,thick] (0.2,0.5) node {$\hat{\Omega}$};
	\draw[black, thick] (0.5, 1.1) node {$\overset{\mathcal{A}_f(t)}{\curvearrowright}$};
\end{scope}
\end{tikzpicture}
\caption{Example: domain reference configuration $\hat{\Omega}$ (left) and domain original configuration at time $t$, $\Omega(t)$ (right). In blue we have the fluid domain, in red the solid domain.}
\label{domain deformation}
\end{figure}
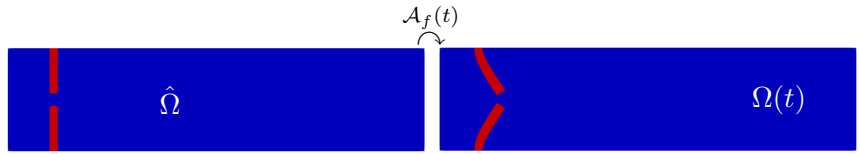
In Figure \ref{domain deformation} we can see an example of a reference configuration and the configuration of the domain at the current time $t$.
Let us see in the next paragraph how to introduce the ALE formalism; for a more detailed discussion about different approaches to describe coupled systems we refer to \cite{Gurtin, Richter}.\\
Let $[0, T]$ be a time interval, and let $\hat{\Omega}_f$ be a \emph{reference configuration} for the fluid.
\begin{definition}
The ALE mapping $\mathcal{A}_f(t)$, for every $t\in[0, T]$ is defined as follows:
\begin{equation*}
\begin{split}
\mathcal{A}_f(t)\colon\hat{\Omega}_f&\mapsto\Omega_f(t)\\
\hat{\bm{x}} &\mapsto \bm{x} = \hat{\bm{x}} + \hat{\bm{d}}_f(\hat{\bm{x}}, t), 
\end{split}
\end{equation*}
\end{definition}
where $\hat{\bm{d}}_f(t)\colon\hat{\Omega}_f\mapsto\mathbb{R}^2$ is the \emph{mesh displacement}. There are different possibilities for the definition of the mesh displacement: in this manuscript, we decide to define $\hat{\bm{d}}_f$ as an harmonic extension of the solid displacement $\hat{\bm{d}}_s$ on the whole fluid domain $\hat{\Omega}_f$: 
\begin{equation*}
\begin{cases}
-\hat{\Delta}\hat{\bm{d}}_f = 0 \quad \text{in }\hat{\Omega}_f,\\
\hat{\bm{d}}_f = \hat{\bm{d}}_s \quad \text{on }\hat{\Gamma}_{FSI},
\end{cases} 
\end{equation*}
{and homogeneous Dirichlet boundary conditions on the remaining portion of the boundary.} Here $\hat{\Gamma}_{FSI}$ is the fluid--structure interface in the reference configuration. A great deal of attention has to be paid to the definition of the mesh displacement, as different choices for $\hat{\bm{d}}_f$ lead to different degrees of regularity: if we lose regularity due to the mesh displacement, we consequently lose regularity at the FSI interface, which is exactly where the coupling between the two physics takes place. It is beyond the scope of this work to discuss the regularity of different definitions of the mesh displacement; nonetheless we refer the interested reader to Chapter $5.3.5$ of \cite{Richter}.
\begin{remark}
$\hat{\bm{d}}_f$ represents the displacement of the grid points, therefore it is not a quantity with a real physical meaning, but rather a geometrical quantity that describes the deformation of the mesh, according to the deformation of the physical domain. It is also important to underline that $\partial_t\hat{\bm{d}}_f\neq \hat{\bm{u}}_f$: in fact, while $\hat{\bm{u}}_f$ represents the velocity of the fluid, $\partial_t\hat{\bm{d}}_f$ is again a geometrical quantity, that can be interpreted as the velocity with which the mesh moves.
\end{remark}

Let us now define the gradient $\bm{F}$ of the ALE map and its determinant $J$, respectively:
\begin{equation*}
\bm{F}:=\hat{\nabla}\mathcal{A}_f, \qquad J:=\text{det}\bm{F}.
\end{equation*}
With these quantities we are ready to present the strong formulation of the FSI problem of interest, within an ALE formalism.
\section{First test case: time dependent FSI problem}\label{problem formulation}
We present the first FSI problem of interest: a time--dependent, non parametrized, nonlinear multiphysics test case.
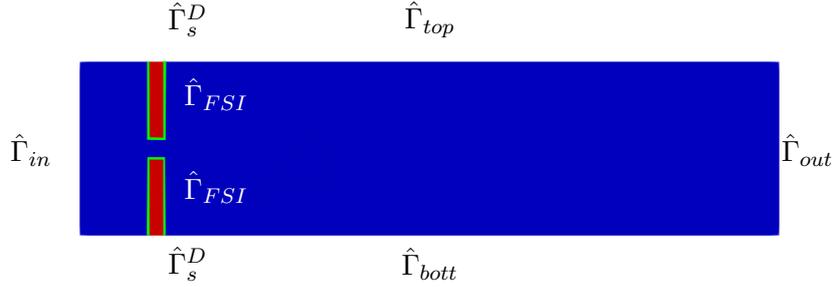
\begin{figure}
\centering
\begin{tikzpicture}
\node[anchor=south west,inner sep=0] (image) at (0,0) {\includegraphics[width=0.6\textwidth]{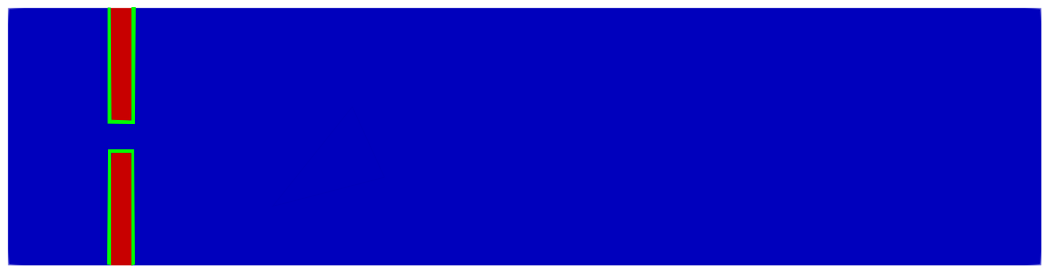}};
\begin{scope}[x={(image.south east)},y={(image.north west)}]
	\draw[black,thick] (-0.06,0.5) node {$\hat{\Gamma}_{in}$};
	\draw[black,thick] (1.03,0.5) node {$\hat{\Gamma}_{out}$};
	\draw[white,thick] (0.2,0.8) node {$\hat{\Gamma}_{FSI}$};
	\draw[white,thick] (0.2,0.3) node {$\hat{\Gamma}_{FSI}$};
	\draw[black,thick] (0.16,1.2) node {$\hat{\Gamma}_s^D$};
	\draw[black,thick] (0.16,-0.11) node {$\hat{\Gamma}_s^D$};
	\draw[black,thick] (0.5, 1.2) node {$\hat{\Gamma}_{top}$};
	\draw[black,thick] (0.5, -0.11) node {$\hat{\Gamma}_{bott}$};
\end{scope}
\end{tikzpicture}
\caption{Physical reference configuration. Blue domain: the reference fluid configuration $\hat{\Omega}_f$. Red leaflets: the reference solid configuration $\hat{\Omega}_s$. The fluid--structure interface $\hat{\Gamma}_{FSI}$ is depicted in green. $\hat{\Gamma}_s^D$: the part of the leaflets that does not move. {$\hat{\Gamma}_{top}$ and $\hat{\Gamma}_{bott}$ are the fluid channel top and bottom walls, $\hat{\Gamma}_{in}$ is the fluid inlet boundary, $\hat{\Gamma}_{out}$ is the fluid outlet boundary.}}
\label{leaflets_ALE_domain}
\end{figure}
{The goal is to simulate the behavior of an incompressible, Newtonian fluid interacting with a linear elastic solid, in the time interval $[0, T]$;} Figure \ref{leaflets_ALE_domain} shows the physical domain in its reference configuration. 
\subsection{Strong formulation}
The coupled FSI problem, formulated over the original configuration, reads as follows: 
find $\bm{u}_f\colon\Omega_f(t)\mapsto\mathbb{R}^2$, $p_f\colon\Omega_f(t)\mapsto\mathbb{R}$ and $\hat{\bm{d}}_s\colon\hat{\Omega}_s\mapsto\mathbb{R}^2$ such that:
\begin{equation}
\label{strong formulation in current configuration}
\begin{cases}
\rho_f(\partial_t\bm{u}_f + (\bm{u}_f\cdot{\nabla})\bm{u}_f)  - \text{div}\sigma_f(\bm{u}_f, p_f) = 0 \quad\text{in }\Omega_f(t)\times(0, T],\\
\text{div}\bm{u}_f = 0 \quad \text{in }\Omega_f(t)\times(0, T],\\
\rho_s\partial_{tt}\hat{\bm{d}_s} -\hat{\text{div}}(\hat{\bm{P}}(\hat{\bm{d}_s})) = 0\quad \text{in }\hat{\Omega}_s\times(0, T],
\end{cases}
\end{equation}
In system $\eqref{strong formulation in current configuration}$, the $\hat{\text{div}}$ denotes the fact that the divergence is computed with respect to $\hat{\bm{x}}$, the space variable in the reference configuration. $\rho_f$ and $\rho_s$ are the fluid and the solid density, while $\sigma_f$ is the fluid \emph{Cauchy stress tensor} for an incompressible Newtonian fluid, and 
$\hat{\bm{P}}$ is the solid \emph{first Piola--Kirchoff tensor} (compressible, linear elastic solid). 
{System \eqref{strong formulation in current configuration} is then completed by some initial conditions (we assume the system to be at rest at the starting time of the simulation), and by the following boundary conditions:
\begin{equation*}
\begin{cases}
{\sigma}_f({\bm{u}}_f, {p}_f){\bm{n}} = -p_{in}(t)\hat{\bm{n}} \text{ on ${\Gamma}_{in}$},\\
{\sigma}_f({\bm{u}}_f, {p}_f){\bm{n}} = 0\text{ on ${\Gamma}_{out}$},\\
\hat{\bm{d}}_s = 0\quad\text{on $\hat{\Gamma}_s^D$},\\
{{\bm{u}}_f=0 \text{ on ${\Gamma}_{top}\cup \Gamma_{bott}$}},
\end{cases}
\end{equation*}
where $p_{in}(t)$ is a time--dependent pressure pulse, $\bm{n}$ is the outward unit normal to the boundary being considered,  ${\Gamma}_{in}$ and $\Gamma_{out}$ represent the inlet and outlet boundaries depicted in Figure \ref{leaflets_ALE_domain}, and $\hat{\Gamma}_s^D$ is the portion of the leaflets' boundary that is attached to the top and bottom walls $\Gamma_{top}$ and $\Gamma_{bott}$. 
It only remains to state the coupling conditions that take place at the fluid--structure interface:
\begin{equation}
\label{coupling conditions current domain}
\begin{cases}
\bm{u}_f = \frac{d}{dt}\bm{d}_s \text{ on }\Gamma_{FSI}(t),\\
\sigma_f(\bm{u}_f, p_f)\bm{n}_f = - J_s^{-1}\hat{\bm{P}}\bm{F}_s^{T}\bm{n}_s \text{ on }\Gamma_{FSI}(t),
\end{cases}
\end{equation}
where $\bm{d}_s$ is the solid deformation in the current configuration $\Omega_s(t)$, $\bm{F}_s:=\hat{\nabla}\hat{\bm{d}_s} + \bm{I}$, $\bm{I}$ the ${2\times 2}$ identity matrix, $J_s:=\text{det}\bm{F}_s$, $\bm{n}_f$ and $\bm{n}_s$ are the unit normals to the FSI interface $\Gamma_{FSI}(t)$ outgoing the fluid and the solid domain, respectively. In system \eqref{coupling conditions current domain} we can interpret the first condition as a {kinematic coupling}, which requires the continuity of the velocities at the interface (the fluid sticks to the moving boundary), whereas the second equation is a {dynamic coupling} corresponding to an action--reaction principle, which is simply stating that the two stresses have to balance out at the interface.}
\\With the formalism introduced in Section \ref{ALE formalism}, we are able to perform a pull--back of the fluid equations onto the fluid reference configuration $\hat{\Omega}_f$; 
the FSI problem on the reference configuration $\hat{\Omega}:=\hat{\Omega}_f\cup\hat{\Omega}_s$ now reads: for every $t\in [0, T]$, find the fluid velocity $\hat{\bm{u}}_f(t):\hat{\Omega}_f\mapsto\mathbb{R}^2$, the fluid pressure $\hat{p}_f(t):\hat{\Omega}_f\mapsto\mathbb{R}$, the fluid displacement $\hat{\bm{d}}_f(t)\colon\hat{\Omega}_f\mapsto\mathbb{R}^2$ and the solid deformation $\hat{\bm{d}}_s(t):\hat{\Omega}_s\mapsto\mathbb{R}^2$ such that
\begin{equation}
\begin{cases}
\rho_fJ(\partial_t\hat{\bm{u}}_f + \hat{\nabla}\hat{\bm{u}}_f\bm{F}^{-1}(\hat{\bm{u}}_f-\partial_t\hat{\bm{d}}_f))  - \hat{\text{div}}(J\hat{\sigma}_f(\hat{\bm{u}}_f, \hat{p}_f)\bm{F}^{-T}) =0 \quad\text{in }\hat{\Omega}_f\times(0, T],\\
\hat{\text{div}}(J\bm{F}^{-1}\hat{\bm{u}}_f) = 0 \quad \text{in }\hat{\Omega}_f\times(0, T],\\
-\hat{\Delta}\hat{\bm{d}}_f = 0 \quad \text{in }\hat{\Omega}_f\times(0, T],\\
\rho_s\partial_{tt}\hat{\bm{d}_s} -\hat{\text{div}}(\hat{\bm{P}}(\hat{\bm{d}_s})) =0\quad \text{in }\hat{\Omega}_s\times(0, T].
\end{cases}
\label{non parametric FSI}
\end{equation}
The fluid tensor $\hat{\sigma}_f$ is the representation in the reference configuration of the Cauchy stress tensor:
\begin{equation*}
\hat{\sigma}_f(\hat{\bm{u}}_f, \hat{p}_f) = \mu_f(\hat{\nabla}\hat{\bm{u}}_f\bm{F}^{-1} + \bm{F}^{-T}\hat{\nabla}^T\hat{\bm{u}}_f) -\hat{p}_f\bm{I},
\end{equation*}
and 
\begin{equation*}
\begin{split}
\hat{\bm{P}}(\hat{\bm{d}}_s) &= \lambda_s\text{tr}\varepsilon_s(\hat{\bm{d}}_s)\bm{I} + 2\mu_s\varepsilon_s(\hat{\bm{d}}_s),\\
\varepsilon_s(\hat{\bm{d}}s) &= \frac{1}{2}(\hat{\nabla}\hat{\bm{d}}_s + \hat{\nabla}^T\hat{\bm{d}}_s).
\end{split}
\end{equation*}
where $\lambda_s$ and $\mu_s$ are the first and second Lam\'e constant of the solid, respectively ($\mu_s$ is also referred to as \emph{shear modulus}).
System \eqref{non parametric FSI} is completed by the same initial conditions, by the boundary conditions
\begin{equation}
\label{boundary conditions}
\begin{cases}
J\hat{\sigma}_f(\hat{\bm{u}}_f, \hat{p}_f)\bm{F}^{-T}\hat{\bm{n}} = -p_{in}(t)\hat{\bm{n}}\quad\text{on $\hat{\Gamma}_{in}$},\\
J\hat{\sigma}_f(\hat{\bm{u}}_f, \hat{p}_f)\bm{F}^{-T}\hat{\bm{n}} = 0\quad\text{on $\hat{\Gamma}_{out}$},\\
\hat{\bm{d}}_s = 0\quad\text{on $\hat{\Gamma}_s^D$},\\
{\hat{\bm{u}}_f=0 \text{ on $\hat{\Gamma}_{top}\cup \hat{\Gamma}_{bott}$}}
\end{cases}
\end{equation}
and by the following coupling conditions:
\begin{equation}
\label{FSI non parametric: coupling conditions}
\begin{cases}
\hat{\bm{d}}_f = \hat{\bm{d}}_s \quad\text{on }\hat{\Gamma}_{FSI}\\ 
\hat{\bm{u}}_f = \partial_t\hat{\bm{d}}_s\quad \text{on }\hat{\Gamma}_{FSI},\\
J\hat{\sigma}_f(\hat{\bm{u}}_f, \hat{p}_f)\bm{F}^{-T}\hat{\bm{n}}_f = -\hat{\bm{P}}(\hat{\bm{d}}_s)\hat{\bm{n}}_s \quad\text{on }\hat{\Gamma}_{FSI}.
\end{cases}
\end{equation}
In the previous equations the vector $\hat{\bm{n}}$ represents the normal vector to the inlet (or outlet) boundary in the reference configuration, whereas $\hat{\bm{n}}_f$ and $\hat{\bm{n}}_s$ are the unit normals to the FSI interface $\hat{\Gamma}_{FSI}$, outgoing the fluid and the solid domain, respectively. 
{In system \eqref{FSI non parametric: coupling conditions}, the first equation is a geometric condition, which states that the fluid and the solid domain do not overlap. The second condition is the {kinematic condition}, expressed in the reference configuration, and similarly the third condition is the {dynamic condition} in the reference configuration.}
\begin{remark}
The gradient and the divergence in equation \eqref{non parametric FSI} are computed with respect to the spatial coordinates in the reference configuration, namely $\hat{\bm{x}}$. Nevertheless, from now on, since everything will be formulated and computed on the reference configuration, in order to ease the exposition, we will drop the $\hat{}$ notation. 
\end{remark}

\subsection{Offline computational phase}\label{offline phase}
We are now going to describe the offline phase of the partitioned procedure that we use to solve the FSI problem of this section. The algorithm is based on a Chorin-Temam projection scheme for the incompressible Navier--Stokes equations \cite{GuermondQuartapelle, Guermond}, and we choose to treat the coupling conditions \eqref{FSI non parametric: coupling conditions} in a semi--implicit way (see also \cite{Ballarin2017, Quaini2008, FernandezGerbeauGrandmont}). We first apply a time stepping procedure to design the algorithm, and then we show the space discretization of the whole procedure.
\subsubsection{High fidelity semi--implicit scheme}
We present the offline phase of the partitioned procedure: we use an operator splitting approach, based on a Chorin-Temam projection scheme with pressure Poisson formulation. 
Let $\Delta T$ be a time--step: we discretize the time interval $[0, T]$ with an equispaced sampling $\{t_0, \hdots, t_{N_T}\}$, where $t_i = i\Delta T$, for $i=0, \hdots, N_T$ and $N_T=\frac{T}{\Delta T}$. We discretize the partial derivative of a function $f$ with a first backward difference BDF1:
\begin{equation*}
D_tf^{i+1} = \frac{f^{i+1} - f^i}{\Delta T}, \quad D_{tt}f^{i+1}= D_t(D_tf^{i+1}),
\end{equation*}
where $f^{i+1}=f(t^{i+1})$. {Hereafter, we will make use of the BDF1 time discretization for both the fluid and the solid problem. We consider the following semi--implicit time discretization of \eqref{non parametric FSI}: for $i = 0, \dots, N_T-1$}
\subsubsection*{\underline{Extrapolation of the mesh displacement $\bm{d}_f$}:} 
find $\bm{d}_f^{i+1}\colon\Omega_f\mapsto\mathbb{R}^2$ such that:
\begin{equation}
\label{mesh displacement}
\begin{cases}
-\Delta\bm{d}_f^{i+1}=0 \quad\text{in $\Omega_f$,}\\
\bm{d}_f^{i+1} = \bm{d}_s^i \quad\text{on $\Gamma_{FSI}$}. 
\end{cases}
\end{equation}
{In this step we are imposing the first of our three coupling conditions, namely the continuity of the displacements at the fluid--structure interface.  This condition is imposed in an explicit way, with respect to time, because we are taking into account the solid displacement $\bm{d}_s^i$ at the previous time iteration $i$, and the mesh displacement unknown $\bm{d}_f^{i+1}$ at the current time iteration $i+1$.}
{The choice of imposing the continuity of the displacements in an explicit way is inspired by many works present in the literature of FSI, see for example \cite{Ballarin2017}.}
{
\begin{remark}
Here we treat the first coupling condition in an explicit way: this approach is less strong than a monolithic one, where this coupling condition would be imposed weakly for the fluid and solid displacement at the same timestep $t^{i+1}$ (see for example \cite{monica}). Nevertheless this choice will allow us to build the fluid displacement $\bm{d}_f$ in a cheap way in the online phase.
\end{remark}
}
\subsubsection*{\underline{Fluid explicit step}:}
find $\bm{u}_f^{i+1}\colon\Omega_f\mapsto\mathbb{R}^2$ such that:
\begin{equation}
\label{label}
\begin{cases}
\begin{split}
&J\rho_f \Bigl(D_t\bm{u}_f^{i+1} + \nabla\bm{u}_f^{i+1}\bm{F}^{-1}(\bm{u}_f^{i+1} -D_t\bm{d}_f^{i+1})\Bigr) -\mu_f\text{div}(J\varepsilon(\bm{u}_f^{i+1})\bm{F}^{-T}) + \\
& + J\bm{F}^{-T}\nabla p_f^i= 0 \text{ in $\Omega_f$,}\\
\end{split}
\\
\bm{u}_f^{i+1} = D_t\bm{d}_f^{i+1} \quad \text{on $\Gamma_{FSI}$,}
\end{cases}
\end{equation}
with $\varepsilon(\bm{u}_f^{i+1}):= \mu_f({\nabla}{\bm{u}}_f\bm{F}^{-1} + \bm{F}^{-T}{\nabla}^T{\bm{u}}_f)$. {Here, we are imposing the dynamic condition (continuity of the velocities at the FSI interface), again in an explicit fashion with respect to time, since now the fluid displacement $\bm{d}_f^{i+1}$ is already known.}
\subsubsection*{\underline{Implicit step}:}
\begin{itemize}
\item \textbf{fluid projection substep (pressure Poisson formulation):} find $p_f^{i+1}\colon\Omega_f\mapsto\mathbb{R}^2$ such that:
\begin{equation}
\label{pressure implicit step}
\begin{cases}
-\text{div}(J\bm{F}^{-1}\bm{F}^{-T}\nabla p_f^{i+1})=-\frac{\rho_f}{\Delta t}\text{div}(J\bm{F}^{-1}\bm{u}_f^{i+1})  \quad\text{in $\Omega_f$,}\\
- \bm{F}^{-T}\nabla p_f^{i+1}\cdot J\bm{F}^{-T}\bm{n}_f = \rho_fD_{tt}\bm{d}_s^{i+1}\cdot J\bm{F}^{-T}\bm{n}_f \quad \text{on $\Gamma_{FSI}$,}
\end{cases}
\end{equation}
subject to the boundary conditions:
\begin{equation}
\label{pressure BC}
\begin{cases}
p_f^{i+1} = p_{in}(t^{i+1}) \quad\text{on }\Gamma_{in} \\
p_f^{i+1} = 0 \quad\text{on }\Gamma_{out},
\end{cases}
\end{equation}
\item \textbf{structure projection substep:} find $\bm{d}_s^{i+1}\colon\Omega_s\mapsto\mathbb{R}^2$ such that:
\begin{equation}
\label{solid implicit step}
\begin{cases}
J_s\rho_sD_{tt}\bm{d}_s^{i+1} - \text{div}(\bm{P}(\bm{d}_s^{i+1}))=0 \quad \text{in }\Omega_s,\\
J\sigma_f(\bm{u}_f^{i+1}, p_f^{i+1})\bm{F}^{-T}\bm{n}_f = -\bm{P}(\bm{d}_s^{i+1})\bm{n}_s\quad \text{on } \Gamma_{FSI}.
\end{cases}
\end{equation}
\end{itemize}
subject to the boundary condition $\bm{d}_s^{i+1}=0$ on $\Gamma_s^D$.
\begin{remark}
We remark that the time--stepping schemes for the fluid problem and for the solid problem are implicit. The denomination "semi--implicit" comes from the fact that the coupling conditions are treated differently. Indeed in system \eqref{label}, the geometrical coupling condition (the second equation) is treated explicitly; on the other hand, the coupling on the fluid and solid velocity (second equation in system\eqref{pressure implicit step}), as well as the coupling of the stresses at the fluid--structure interface (second equation in system \eqref{solid implicit step}), are treated implicitly.
\end{remark}
The implicit step couples pressure stresses to the structure, and it is iterated until convergence is reached. {We would like to remark that, throughout this manuscript, the BDF1 time stepping scheme is used also for the solid problem: there are other alternatives, such as a Newmark--Beta or HHT-alpha methods, that represent a standard choice in solid mechanics, as they have some desirable properties regarding stability and dissipation. In our work we did not encounter any problem with the stability in time of the algorithm, and therefore chose to use a BDF1 scheme for its easiness of implementation.}
\begin{remark}
In the implicit step \eqref{pressure implicit step} we have chosen a pressure Poisson formulation; an alternative is to use a Darcy formulation, which is defined as follows:
find $p_f^{i+1}$ and $\tilde{\bm{u}}_f^{i+1}$ such that:
\begin{equation*}
\begin{cases}
\rho_fJ\frac{\tilde{\bm{u}}_f^{i+1}-\bm{u}_f^{i+1}}{\Delta T} + J\bm{F}^{-T}\nabla p_f^{i+1} = 0\quad \text{in }\Omega_f,\\
\text{div}(J\bm{F}^{-1}\tilde{\bm{u}}_f^{i+1}) = 0 \quad \text{in $\Omega_f$}.
\end{cases}
\end{equation*}
Throughout this manuscript we choose to employ a Poisson formulation, for the sake of a more efficient reduced order model, since the Darcy formulation requires the introduction of an additional unknown $\tilde{\bm{u}}_f$, which translates in a larger system, comprised of both velocity and pressure, at the implicit step.
\end{remark}
{Let us now have a look at  equations \eqref{label}--\eqref{solid implicit step}: the fluid problem is solved using Dirichlet boundary conditions (the displacement computed at the previous timestep), and the solid problem is solved using Neumann boundary conditions (the fluid normal stress just computed). However, as it is mentioned in \cite{BadiaNobileVergara} and references therein, these kind of partitioned schemes (Dirichlet--Neumann couplings) usually require a large amount of sub--iterations of the implicit step, before a convergence between fluid and solid problem is reached, especially in those situations where the added mass effect is particularly heavy (e.g. blood flow simulations). Motivated by this, in order to have a better control on the number of sub--iterations needed in our algorithm, we decide to replace  the Neumann condition is system \eqref{pressure implicit step} with a Robin coupling condition, as suggested by  \cite{AstorinoChoulyFernandez, Ballarin2017, BadiaNobileVergara}. In \cite{Ballarin2017} the authors propose a Robin coupling condition that is based on a coefficient $\alpha_{ROB}$ that has been computed for the one dimensional structure. For our problem however, the solid is two dimensional and elastic: we therefore rely on the work presented in \cite{BanksHenshaw}, where the authors compute the constant $\alpha_{ROB}$ in the case of an elastic solid. We therefore just have to incorporate the expression of $\alpha_{ROB}$ found in \cite{BanksHenshaw} into the Robin coupling condition presented in \cite{Ballarin2017}, and remember to pull back the condition onto the reference fluid--structure interface, using the ALE map. The final expression of the Robin coupling condition is:
}
\begin{equation}
\label{Robin coupling}
\alpha_{ROB}p^{i+1} + \bm{F}^{-T}\nabla p^{i+1}\cdot J\bm{F}^{-T}\bm{n}_f = \alpha_{ROB} p^{i+1,\star} -\rho_fD_{tt}\bm{d}_s^{i+1,\star}\cdot J\bm{F}^{-T}\bm{n}_f.
\end{equation}
In equation \eqref{Robin coupling}, $p^{i+1, \star}$ and $\bm{d}_s^{i+1, \star}$ are suitable extrapolations of the fluid pressure and the solid displacement, respectively; we show in the next paragraph which kind of extrapolation we use. The constant $\alpha_{ROB}$ is defined as $\alpha_{ROB} = \frac{\rho_f}{z_p\Delta T}$ where $z_p$ is called the \emph{solid impedance}:
\begin{equation*}
\begin{split}
z_p &= \rho_s c_p,\\
c_p &= \sqrt{\frac{\lambda_s + 2\mu_s}{\rho_s}}.
\end{split}
\end{equation*}
{Condition \eqref{Robin coupling} is imposed only on the fluid side: Robin conditions are indeed nonstandard in solid mechanics, therefore a lot of already existing codes would not allow to impose such a condition for the solid problem.}
\subsubsection{Space discretization of the semi--implicit procedure}\label{space discretization section}
We now present the semi--discretized version of the algorithm introduced. We define the following function spaces for the fluid:
\begin{equation*}
V(\Omega_f) := [H^1(\Omega_f)]^2, \quad E^f(\Omega_f) := [H^1(\Omega_f)]^2, \quad Q(\Omega_f):= L^2(\Omega_f),
\end{equation*}
endowed with the $H^1$ norm ($V(\Omega_f)$), the $H^1$ seminorm ($E^f(\Omega_f)$)and the $L^2$ norm respectively, and the function space for the solid: $E^s(\Omega_s)= [H^1(\Omega_s)]^2$, endowed with the $H^1_0$ norm. {We discretize in space the FSI problem, using second order Lagrange Finite Elements for the fluid velocity, resulting in the discrete space $V_h(\Omega_f)\subset V(\Omega_f)$, while for the fluid pressure, the fluid displacement and the solid displacement we use first order Lagrange Finite Elements, resulting in the discrete space $Q_h(\Omega_f)\subset Q(\Omega_f)$, $E^f_h(\Omega_f)\subset E^f(\Omega_f)$ and $E_h^s(\Omega_s)\subset E^s(\Omega_s)$; we further assume here that the fluid and the solid discretizations match at the FSI interface.}
The non--homogeneous boundary condition in system \eqref{pressure BC} can be easily treated by introducing, at timestep $t^{i+1}$, a lifting function $\ell^{i+1}$ such that $\ell^{i+1}= p_{in}(t^{i+1})$ on $\Gamma_{in}$ and $\ell^{i+1}=0$ on $\Gamma_{out}$; we refer, for example, to \cite{Supremizer, monica} for more details concerning the use of a lifting function within model order reduction. By introducing the homogenized pressure $p_f^{0, i+1}:=p_f^{0, i+1} - \ell^{i+1}$, we can conclude now that $p_f^{0, i+1}\in Q_h^0$, where $Q_h^0 = \{q_h\in Q_h\colon q_h=0\quad\text{on }\Gamma_{in}\cup\Gamma_{out}\}$.
The discretized version of the semi--implicit procedure reads as follows: for $i=0, \hdots, N_T$,
\subsubsection*{\underline{Extrapolation of the mesh displacement}:} 
find $\bm{d}_{f, h}^{i+1}\in E_h^f$ such that $\forall \bm{e}_{f,h} \in E^f_h$:
\begin{equation}
\label{mesh displacement}
\begin{cases}
\int_{\Omega_f}\nabla\bm{d}_{f, h}^{i+1}\cdot\nabla\bm{e}_{f, h}\,dx =0,\\
\bm{d}_{f, h}^{i+1} = \bm{d}_{s, h}^i \quad\text{on $\Gamma_{FSI}$}. 
\end{cases}
\end{equation}
\subsubsection*{\underline{Fluid explicit step}:}
find $\bm{u}_{f, h}^{i+1}\in V_h$ such that $\forall \bm{v}_h\in V_h$:
\begin{equation}
\label{fluid explicit step test1}
\begin{cases}
\begin{split}
&\rho_f\int_{\Omega_f} JD_t\bm{u}_{f, h}^{i+1}\cdot\bm{v}_h\,dx 
+ \rho_f\int_{\Omega_f}J(\nabla\bm{u}_{f, h}^{i+1}\bm{F}^{-1}(\bm{u}_{f, h}^{i+1} -D_t\bm{d}_{f, h}^{i+1}))\cdot\bm{v}_h\,dx\\
&+ \mu_f\int_{\Omega_f}J\varepsilon(\bm{u}_{f, h}^{i+1})\bm{F}^{-T}:\nabla\bm{v}_h\,dx + \int_{\Omega_f}J\bm{F}^{-T}\nabla p_{f, h}^i\cdot\bm{v}_h\,dx= 0
\end{split}
\\ \bm{u}_{f, h}^{i+1} = D_t\bm{d}_{f, h}^{i+1} \quad \text{on $\Gamma_{FSI}$,}
\end{cases}
\end{equation}
{This step results in a nonlinear system, which, at the computational level, is being solved with a Newton method. We remark that in the weak formulation, the boundary terms for the velocity vanish: this is in part due to the homogeneous Dirichlet boundary condition on $\Gamma_{top}\cup \Gamma_{bott}$, and in part due to the fact that, in the Chorin--Temam scheme, the inlet boundary condition (first equation in system \eqref{boundary conditions}) is split between the velocity and the pressure in the following way: $J\varepsilon(\bm{u}_{f, h}^{i+1})\bm{F}^{-T}\bm{n}=0$ on $\Gamma_{in}$ and $Jp_{f, h}^{i+1}\bm{F}^{-T}\bm{n}=-p_{in}$ on $\Gamma_{in}$. The imposition of Neumann boundary conditions within a Chorin--Temam scheme is not at all trivial, and we refer the interested reader to a more detailed discussion presented in \cite{Rannacher, KARNIADAKIS1991414}.}
\subsubsection*{\underline{Implicit step}:} for any $j=0, \dots$ until convergence:
\begin{enumerate}
\item \textbf{fluid projection substep (pressure Poisson formulation):} find $p_{f, h}^{0, i+1, j+1}\in Q_h^0$ such that $\forall q_h\in Q_h^0$:
\begin{equation*}
\begin{split}
&\alpha_{ROB}\int_{\Gamma_{FSI}}p_{f, h}^{0, i+1, j+1}q_h\,ds + \int_{\Omega_f}J\bm{F}^{-T}\nabla p_{f, h}^{0, i+1, j+1}\cdot \bm{F}^{-T}\nabla q_h\,dx = \\
&-\frac{\rho_f}{\Delta T}\int_{\Omega_f}\text{div}(JF^{-1}\bm{u}_{f, h}^{i+1})q_h\,dx - \rho_f\int_{\Gamma_{FSI}}D_{tt}\bm{d}_{s, h}^{i+1, j}\cdot J\bm{F}^{-T}\bm{n}_fq_h\,ds\\
&+ \alpha_{ROB}\int_{\Gamma_{FSI}}p_{f, h}^{i+1, j}q_h\,ds  - \alpha_{ROB}\int_{\Gamma_{FSI}}\ell^{i+1}q_h\,ds- \int_{\Omega_f}J\bm{F}^{-T}\nabla \ell^{i+1}\cdot \bm{F}^{-T}\nabla q_h\,dx,\\
\end{split}
\end{equation*}
\item \textbf{structure projection substep:} find $\bm{d}_{s, h}^{i+1, j+1}\in E^s_h$ such that $\forall \bm{e}_{s, h}\in E^s_h$:
\begin{equation*}
\begin{split}
&\rho_s\int_{\Omega_s}J_sD_{tt}\bm{d}_{s, h}^{i+1, j+1}\cdot\bm{e}_{s, h}\,dx + \int_{\Omega_s}\bm{P}(\bm{d}_{s, h}^{i+1, j+1}):\nabla\bm{e}_{s, h}\,dx =-\\
&-\int_{\Gamma_{FSI}}J\sigma_f(\bm{u}_{f, h}^{i+1}, p_{f, h}^{i+1, j+1})\bm{F}^{-T}\bm{n}_f\cdot\bm{e}_{s, h}\,dx,
\end{split}
\end{equation*}
subject to the boundary condition $\bm{d}_{s, h}^{i+1, j+1} = 0$ on $\Gamma_{D}^s$.
\end{enumerate}
We iterate between the two implicit substeps, until a convergence criteria is satisfied; we choose as stopping criteria a relative error on the increments of the pressure and the solid displacement, namely:
\begin{equation*}
\text{max}\Bigl(\frac{\lvert\lvert p_{f, h}^{i+1, j+1} - p_{f, h}^{i+1, j}\rvert\rvert_{Q_h}}{\lvert\lvert p_{f, h}^{i+1, j+1}\rvert\rvert_{Q_h}}; \frac{\lvert\lvert \bm{d}_{s, h}^{i+1, j+1} - \bm{d}_{s, h}^{i+1, j}\rvert\rvert_{E^s_h}}{\lvert\lvert \bm{d}_{s, h}^{i+1, j+1}\rvert\rvert_{E^s_h}}\Bigr) < \varepsilon,
\end{equation*}
where $\varepsilon$ is a fixed tolerance.\\
In the pressure Poisson formulation, to impose the Robin coupling condition, we have chosen the pressure at the previous implicit iteration, namely $p_f^{i+1, j}$, as an extrapolation for the fluid pressure, and the same goes for the extrapolation of the structure displacement. 
\subsubsection{POD and reduced basis generation}
For the generation of the reduced basis for the fluid velocity $\bm{u}_f$ and the fluid displacement $\bm{d}_f$ we pursue here the idea that was first proposed in \cite{Ballarin2017}. For the homogenized fluid pressure $p_f^0$ and for the solid displacement $\bm{d}_s$ we employ a standard POD, giving rise to the reduced spaces $Q_N^0$ and $E_N^s$ respectively, {though the authors would like to mention the fact that, as an alternative to the POD modes for the solid problem, the so--called vibrational modes can be used: these are obtained solving a generalized eigenvalue problem involving the mass and the stiffness matrix of the solid problem. Vibrational modes show very good results, especially for linear problems, and the authors refer the interested reader to \cite{Kalashnikova_ROMcomprohtua, osti_974411}}.
\subsubsection*{\underline{Change of variable for the fluid velocity}}
The main idea here is to introduce a change of variable in the fluid problem, in order to transform the non homogeneous Dirichlet condition at the FSI interface in system \eqref{fluid explicit step test1} into a homogeneous boundary condition. The motivation of this choice is that, to impose the second condition in system \eqref{fluid explicit step test1}, we could use a Lagrange multiplier $\lambda$, thus increasing the dimension of the system to be solved in the online phase. In order to avoid this and in order to design a more efficient reduced method, we choose to transform the non--homogeneous coupling condition into a homogeneous one: {we refer to \cite{GUNZBURGER20071030} for a detailed discussion on the treatment of non--homogeneous Dirichlet boundary conditions within a model order reduction framework}.
We begin by defining a new variable $\bm{z}_{f, h}^{i+1}$:
\begin{equation*}
\bm{z}_{f, h}^{i+1} := \bm{u}_{f, h}^{i+1} - D_t\bm{d}_{f, h}^{i+1}.
\end{equation*} 
With this change of variable, the second equation in \eqref{fluid explicit step test1} is now equivalent to the homogeneous boundary condition for the new variable: 
\begin{equation*}
\bm{z}_{f, h}^{i+1}=0 \quad \text{on $\Gamma_{FSI}$},
\end{equation*}
for which no imposition by means of Lagrange multiplier is needed.
Therefore, during the offline phase of the algorithm, at every iteration $i+1$, after we have computed the velocity $\bm{u}_{f, h}^{i+1}$, we compute the change of variable $\bm{z}_{f, h}^{i+1}$. We then consider the following snapshots matrix:
\begin{equation*}
\bm{\mathcal{S}}_z = [{\bm{z}}_{f,h}^1,\hdots,{\bm{z}}_{f,h}^{N_T}] \in \mathbb{R}^{\mathcal{N}_u^h\times N_T},
\end{equation*}
where $\mathcal{N}_u^h = \text{dim}V_h$. 
We then apply a POD to the snapshots matrix $\bm{\mathcal{S}}_{z}$ and we retain the first $N_z$ POD modes $\Phi_{\bm{z}}^1, \hdots, \Phi_{\bm{z}}^{N_z}$. We therefore have the reduced space:
\begin{equation*}
V^N:=\text{span}\{\Phi_{\bm{z}}^k\}_{k=1}^{N_z},
\end{equation*}
and now it is clear that, since every $\Phi_{\bm{z}}^k$ satisfies the condition $\Phi_{\bm{z}}^k=0$ on $\Gamma_{FSI}$, then also every element of $V_N$ will satisfy the same condition.
\subsubsection*{\underline{Harmonic extension of the fluid displacement}}
In order to generate the reduced basis for the fluid displacement $\bm{d}_f$, we pursue again the idea presented in \cite{Ballarin2017}. We start by generating the snapshots matrix related to the solid displacement:
\begin{equation*}
\bm{\mathcal{S}}_{d_s} = [{\bm{d}}_{s,h}^1, \hdots, {{\bm{d}}}_{s,h}^{N_T}] \in \mathbb{R}^{\mathcal{N}_{\bm{d}_s}^h\times N_T},
\end{equation*}
where $\mathcal{N}_{\bm{d}_s}^h = \text{dim}E_h^s$. 
We then apply a POD to the snapshots matrix and retain the first $N_d$ POD modes $\Phi_{\bm{d}_s}^1, \hdots, \Phi_{\bm{d}_s}^{N_d}$, thus defining the reduced space for the solid problem:
\begin{equation*}
E_N^s:=\text{span}\{\Phi_{\bm{d}_s}^k\}_{k=1}^{N_d}.
\end{equation*}
We then employ an harmonic extension of each one of the reduced basis $\Phi_{\bm{d}_s}^k$ to the fluid domain, thus obtaining the functions $\Phi_{\bm{d}_f}^k$ such that:
\begin{equation*}
\begin{cases}
-\Delta \Phi_{\bm{d}_f}^k = 0 \quad\text{in $\Omega_f$}, \\
\Phi_{\bm{d}_f}^k = \Phi_{\bm{d}_s}^k \quad\text{on $\Gamma_{FSI}$}.
\end{cases}
\end{equation*}
{We impose homogeneous Dirichlet boundary conditions on the remaining part of the boundaries.}
We can then define the reduced space for the fluid displacement:
\begin{equation*}
E_N^f:= \text{span}\{\Phi_{\bm{d}_f}^k\}_{k=1}^{N_d}.
\end{equation*}
The reason for for this choice, instead of employing a standard POD on the set of snapshots of $\bm{d}_f$, is given by the fact that we can avoid the introduction of another Lagrange multiplier to impose the non--homogeneous boundary condition present in system \eqref{mesh displacement}. 
We avoid to solve the reduced system related to \eqref{mesh displacement}: instead of solving an harmonic extension problem at every time--step in the online phase, we solve \emph{once and for all} $N_{\bm{d}}$ harmonic extension problems in the expensive offline phase. Then, during the online phase, the reduced fluid displacement will be computed just as a linear combination of the basis $\Phi_{\bm{d}_f}^i$, with coefficients that are the coefficients of the reduced solid displacement at the previous time--step. We will see in the next section the final formulation of the online phase of the algorithm.
{Before moving on, we summarize the offline computational phase, with the aim of helping the reader to better understand the whole procedure so far. Let $i+1$ be the index of the current time iteration:
\begin{enumerate}
\item compute the snapshot $\bm{d}_{f, h}^{i+1}$, using the previously computed snapshot $\bm{d}_{s, h}^i$;
\item solve the fluid explicit part, and find $\bm{u}_{f, h}^{i+1}$ such that $\bm{u}_{f, h}^{i+1}=D_t\bm{d}_{f, h}^{i+1})$ on the FSI interface;
\item compute the fluid viscous stress $\varepsilon(\bm{u}_{f, h}^{i+1})$;
\item compute $\bm{z}_{f, h}^{i+1}=\bm{u}_{f, h}^{i+1}-D_t\bm{d}_{f, h}^{i+1})$;
\item store the snapshot  $\bm{z}_{f, h}^{i+1}$ in the snapshot matrix $\bm{\mathcal{S}}_z$;
\item iterate until tolerance $\varepsilon$ is reached:
\begin{itemize}
\item solve the pressure Poisson problem, using the solid displacement at the previous subiteration, and find $p_{f, h}^{i+1, j+1}$;
\item solve the solid problem, using the fluid stress tensor $\varepsilon(\bm{u}_{f, h}^{i+1})-p_{f, h}^{i+1, j+1}\bm{I}$, and find $\bm{d}_{s, h}^{i+1, j+1}$;
\end{itemize}
\item store the homogenized snapshot  $\bm{p}_{f, h, 0}^{i+1}:=\bm{p}_{f, h}^{i+1} - \ell^{i+1}$ in the snapshot matrix $\bm{\mathcal{S}}_p$;
\item store the snapshot  $\bm{d}_{s, h}^{i+1}$ in the snapshot matrix $\bm{\mathcal{S}}_{d_s}$;
\item POD compression on $\bm{\mathcal{S}}_z$ $\rightarrow$ $\{\Phi_{\bm{z}}^k\}_{k=1}^{N_z}$;
\item POD compression on $\bm{\mathcal{S}}_p$ $\rightarrow$ $\{\Phi_{p}^k\}_{k=1}^{N_p}$;
\item POD compression on $\bm{\mathcal{S}}_{d_s}$ $\rightarrow$ $\{\Phi_{\bm{d}_s}^k\}_{k=1}^{N_d}$;
\item solve $N_d$ harmonic extension problems and find $\{\Phi_{\bm{d}_f}^k\}_{k=1}^{N_d}$.
\end{enumerate}
}
\subsection{Online computational phase}\label{online phase}
We are now ready to present the online formulation of the partitioned procedure, {which is obtained by means of a Galerkin projection over the reduced spaces $V_N$, $E^f_N$, $Q^f_N$ and $E^s_N$}. For every $i=0, \hdots, N_T$, we introduce the reduced functions $\bm{z}_{f, N}^{i+1}$, $p_{f, N}^{0, i+1}$, $\bm{d}_{s, N}^{i+1}$ of the form:
\begin{gather}
\label{z}
\bm{z}_{f, N}^{i+1} = \sum_{k=1}^{N_{\bm{z}}}\underline{\textbf{z}}_k^{i+1}\Phi_{\bm{z}_f}^k,\\
\label{p}
p_{f, N}^{0, i+1} = \sum_{k=1}^{N_{p}}\underline{p}_k^{0, i+1}\Phi_{p}^k,\\
\label{d}
\bm{d}_{s, N}^{i+1} = \sum_{k=1}^{N_{\bm{d}}}\underline{\textbf{d}}_k^{i+1}\Phi_{\bm{d}_s}^k.
\end{gather}
Then the online phase of the partitioned procedure reads as follows:
\subsubsection*{\underline{Mesh displacement}:}
let $\bm{d}_{f, N}^{i+1}$ be defined by the reduced solid displacement at the previous time--step:
\begin{equation}
\label{d_f}
\bm{d}_{f, N}^{i+1} = \sum_{k=1}^{N_{\bm{d}}}\underline{\textbf{d}}_k^i\Phi_{\bm{d}_f}^k;
\end{equation}
\subsubsection*{\underline{Fluid explicit step (with change of variable)}:}
find $z_{f, N}^{i+1}\in V_N$ such that $\forall \bm{v}_N\in V_N$:
\begin{equation*}
\begin{split}
&\rho_f\int_{\Omega_f} J\Bigl(\frac{\bm{z}_{f, N}^{i+1}-\bm{u}_{f, N}^i}{\Delta T}\Bigr)\cdot\bm{v}_N\,dx 
+ \rho_f\int_{\Omega_f}J(\nabla(\bm{z}_{f, N}^{i+1}+ D_t\bm{d}_{f, N}^{i+1})\bm{F}^{-1}\bm{z}_{f, N}^{i+1})\cdot\bm{v}_N\,dx\\
&+ \mu_f\int_{\Omega_f}J\varepsilon(\bm{z}_{f, N}^{i+1})\bm{F}^{-T}:\nabla\bm{v}_N\,dx + \int_{\Omega_f}J\bm{F}^{-T}\nabla p_{f, N}^i\cdot\bm{v}_h\,dx= \\
&-\rho_f\int_{\Omega_f} J\Bigl(\frac{D_t\bm{d}_{f, N}^{i+1}}{\Delta T}\Bigr)\cdot\bm{v}_N\,dx 
- \mu_f\int_{\Omega_f}J\varepsilon(D_t\bm{d}_{f, N}^{i+1})\bm{F}^{-T}:\nabla\bm{v}_N\,dx,
\end{split}
\end{equation*}
we then define the reduced fluid velocity: $\bm{u}_{f, N}^{i+1} = \bm{z}_{f, N}^{i+1} + D_t\bm{d}_{f, N}^{i+1}$.
\subsubsection*{\underline{Implicit step}:} for any $j=0, \dots$ until convergence:
\begin{enumerate}
\item \textbf{fluid projection substep:} find $p_{f, N}^{0, i+1, j+1}\in Q_N^0$ such that $\forall q_N\in Q_N^0$:
\begin{equation*}
\begin{split}
&\alpha_{ROB}\int_{\Gamma_{FSI}}p_{f, N}^{0, i+1, j+1}q_N\,ds + \int_{\Omega_f}J\bm{F}^{-T}\nabla p_{f, N}^{0, i+1, j+1}\cdot \bm{F}^{-T}\nabla q_N\,dx = \\
&-\frac{\rho_f}{\Delta T}\int_{\Omega_f}\text{div}(JF^{-1}\bm{u}_{f, N}^{i+1})q_N\,dx - \rho_f\int_{\Gamma_{FSI}}(D_{tt}\bm{d}_{s, N}^{i+1, j})\cdot J\bm{F}^{-T}\bm{n}_fq_N\,ds\\
&+ \alpha_{ROB}\int_{\Gamma_{FSI}}p_{f, N}^{i+1, j}q_N\,ds - \alpha_{ROB}\int_{\Gamma_{FSI}}\ell^{i+1}q_N\,ds\\
&- \int_{\Omega_f}J\bm{F}^{-T}\nabla \ell^{i+1}\cdot \bm{F}^{-T}\nabla q_N\,dx\\
\end{split}
\end{equation*}
we then recover the reduced fluid pressure $p_{f, N}^{i+1, j+1} = p_{f, N}^{0, i+1, j+1} + \ell^{i+1}$.
\item \textbf{structure projection substep:} find $\bm{d}_{s, N}^{i+1, j+1}\in E^s_N$ such that $\forall \bm{e}_s\in E^s_N$:
\begin{equation*}
\begin{split}
&\rho_s\int_{\Omega_s}D_{tt}\bm{d}_{s, N}^{i+1, j+1}\cdot\bm{e}_N\,dx + \int_{\Omega_s}\bm{P}(\bm{d}_{s, N}^{i+1, j+1}):\nabla\bm{e}_N\,dx = -\\
&-\int_{\Omega_s}J\sigma_f(\bm{u}_{f, N}^{i+1}, p_{f, N}^{i+1, j+1})\bm{F}^{-T}\bm{n}_f\cdot\bm{e}_N\,dx.
\end{split}
\end{equation*}
\end{enumerate}

\subsection{Numerical results}\label{numerical section non parametrized}
We now present some numerical results obtained with the semi--implicit scheme. The reference physical configuration of the problem of interest is the one represented in Figure \ref{leaflets_ALE_domain}: {the geometrical properties of the domain are reported in Table \ref{table leaflets ALE}}; the leaflets are situated $1$ cm downstream the inlet boundary.
For our simulation we used a time--step $\Delta T = 10^{-4}$, and a final time $T=0.05$ $s$, for a total of $N_T = 500$ iterations.
\begin{table}
\centering
\caption{Values for the implementation of the offline phase.}
\label{table leaflets ALE}
\begin{tabular}{|l c|}
\hline
\textbf{Physical constants} & \textbf{Value}\\
$\rho_f$ & $1$ g/cm$^3$\\ 
$\mu_f$ & $0.035$ Poise\\
$\rho_s$ & $1.1$ g/cm$^3$\\
$\mu_s$ & $10^5$\\
$\lambda_s$ & $8\times 10^5$\\
\hline
\textbf{Geometrical constants} & \textbf{Value}\\
channel length & $10$ cm\\
channel height & $2.5$ cm\\
leaflets' length & $1$ cm\\
leaflets' thickness & $0.2$ cm\\
\hline
\textbf{Discretization details} & \textbf{Value}\\
FE displacement order & $1$\\
FE velocity order & $2$\\
FE pressure order & $1$\\
\hline
\end{tabular}
\end{table}

\begin{figure}
\centering
\includegraphics[scale=0.2]{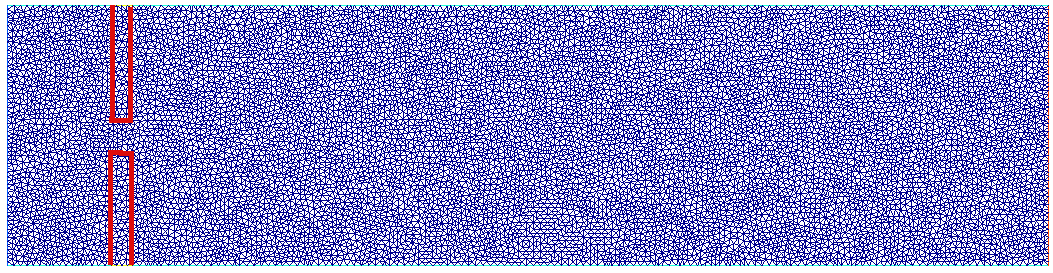}
\caption{{Example of the mesh used for the numerical simulations. The FSI interface has been contoured in red to ease the visualization. {We use a triangular mesh}, where the number of cells is $22060$, for a total of $11269$ vertices.}}
\end{figure}

\begin{figure}[h]
\centering
\begin{minipage}{0.4\textwidth}
\centering
\includegraphics[width=1.\textwidth]{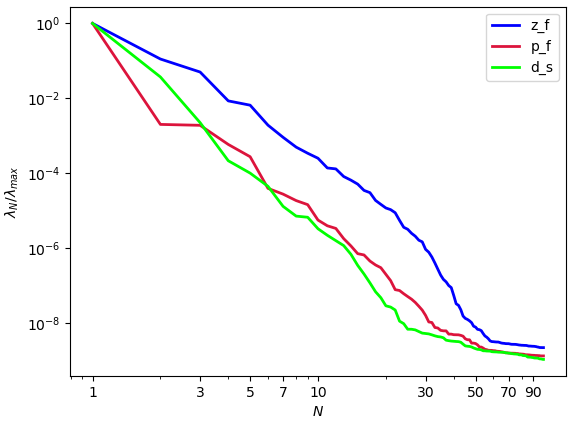}
\subcaption{POD eigenvalues.}\label{leaflets eigenvalues}
\end{minipage}%
\begin{minipage}{0.4\textwidth}
\centering
\includegraphics[width=1.\textwidth]{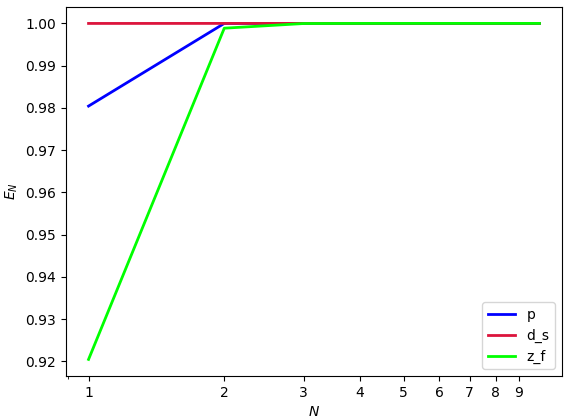}
\subcaption{POD retained energy}\label{retained energy}
\end{minipage}
\caption{Results of the POD for the non parametric problem.}
\end{figure}

\begin{figure}
\centering
\includegraphics[scale=0.3]{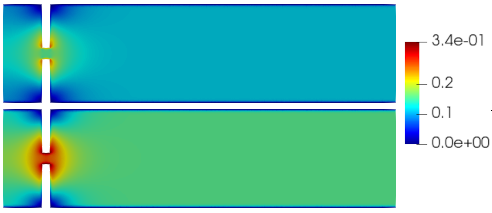}
\caption{Reduced velocity $\bm{u}_{f, N}$ at time--step $t=0.04s$ (top) and at time--step $t=0.05s$ (bottom). The velocity has been obtained with $N_z = 15$ reduced basis.}
\label{reduced velocity non parametric}
\end{figure}

\begin{figure}
\centering
\includegraphics[scale=0.5]{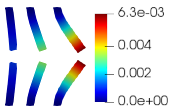}
\caption{Reduced solid displacement $\bm{d}_{s, N}$ at time--step $t=0.035s$ (left), $t=0.04s$ (center) and at time--step $t=0.05s$ (right). The displacement has been obtained with $N_d = 10$ reduced basis. The displacement has been magnified for visualization purposes.}
\label{reduced displacement}
\end{figure}

\begin{figure}
\centering
\includegraphics[scale=0.3]{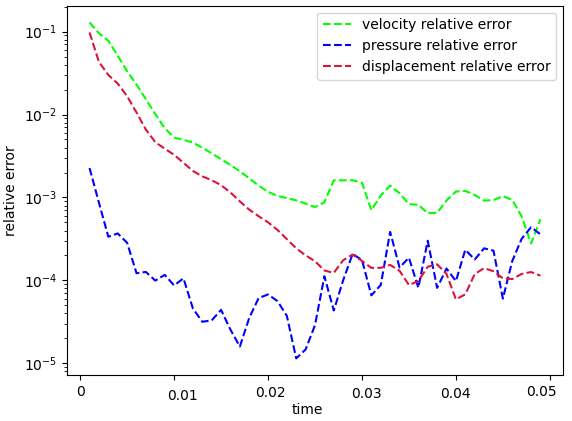}
\caption{Error analysis: relative error behavior, as a function of time. The reduced solutions have been obtained with: $N_{\bm{z}} = 15$, $N_p = 10$ and $N_{\bm{d}_s}=10$.}
\label{error behavior}
\end{figure}

\begin{figure}
\centering
\includegraphics[scale=0.3]{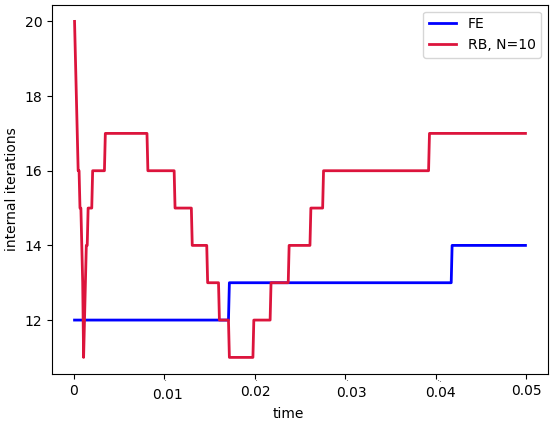}
\caption{Internal iterations: comparison between full order model (blue) and reduced order model (red). The reduced solutions have been obtained with: $N_{\bm{z}} = 15$, $N_p = 10$ and $N_{\bm{d}_s}=10$.}
\label{internal iterations behavior}
\end{figure}

\begin{table}
\centering
\caption{Average relative error for $\bm{u}_f$, $p_f$, $\bm{d}_s$ with basis refinement for the pressure.}
\label{table average relative error non parametric}
\begin{tabular}{|l | c| c| c|}
\hline
$N_p$ & $\bm{u}_f$ & $p_f$ & $\bm{d}_s$\\
\hline
$5$ & $0.010285$ & $0.000199$ & $0.005255$\\ 
$10$ & $0.004370$ & $7.16739\times 10^{-5}$ & $0.002138$\\
$15$ & $0.003013$ & $8.5678\times 10^{-5}$ & $0.002063$\\
{$25$} & {$0.003020$} & {$0.000107$} & {$0.002107$}\\
{$35$} & {$0.0031058$} & {$0.000149$} & {$0.002130$}\\
$40$ & $0.003144$ & $0.000151$ & $0.002078$\\
\hline
\end{tabular}
\end{table}

\begin{table}
\centering
\caption{Average number of subiterations: comparison between full order model, and reduced order model with basis refinement for the pressure.}
\label{table average internal iterations non parametric}
\begin{tabular}{|l | c| c|}
\hline
$N_p$ & full order model & reduced order model\\
\hline
$5$ & $13$ & $15.33$\\ 
$10$ & $13$ & $16$\\
$15$ & $13$ & $26$\\
{$25$} & {$13$} & {$26.37$}\\
{$35$} & {$13$} & {$26.98$}\\
$40$ & $13$ & $26.9$\\
\hline
\end{tabular}
\end{table}

The values of the physical constants used in the simulation are reported in Table \ref{table leaflets ALE}. A pressure impulse $p_{in}(t)$ is applied at the inlet boundary, and after some time this impulse becomes constant:
\begin{equation*}
p_{in}(t)=
\begin{cases}
5 -5\text{cos}\Bigl(\frac{2\pi t}{T_{in}}\Bigr) \qquad \text{for }t\leq 0.025s,\\
5 \qquad \text{for }t>0.025s,
\end{cases}
\end{equation*}
where $T_{in}=0.1$ $s$. We fix a tolerance of $\varepsilon=10^{-6}$ as a stopping criterion for the subiterations between the pressure Poisson problem and the solid problem.

Since we do not consider the top and the bottom walls of the fluid domain to be deformable, we impose a homogeneous boundary condition for the fluid velocity on these walls.

Figure \ref{leaflets eigenvalues} shows the rate of decay of the first $100$ eigenvalues associated with three unknowns of the problem, namely the change of variable for the fluid velocity change of variable $\bm{z}_f$, the pressure $p_f$ and the solid displacement $\bm{d}_s$. It can be noticed that the rate of decay of the eigenvalues for the pressure and for the fluid change of variable is slower than the rate of decay of the eigenvalues of the solid displacement. Moreover, in Figure \ref{retained energy} we can notice that the first mode of the solid displacement retains $2\%$ more energy compared to the first mode of the pressure, and $8\%$ more energy with respect to the first mode of $\bm{z}_f$, which is the one that retains less energy. {Figure \ref{leaflets eigenvalues} is also important to have a first insight on the dimension of the reduced spaces that we are going to take during the online phase: indeed the rate of decay of the eigenvalues returned by the POD gives us an idea of the behavior of the approximation error that we commit by approximating the FE solution with the RB one. The following relationship holds true (we state it for the fluid pressure, but the same holds also for the other components):
\begin{equation*}
{\sum_{i=1}^{N_T}\lvert\lvert p_{f, h}^i - \Pi_{N_p}p_{f, h}^i\rvert\rvert_{L^2(\Omega_f)}^2 = \sum_{k>N_p}\lambda^k_p,}
\end{equation*}
where $\Pi_{N_p}$ is the orthogonal projector onto the POD space of dimension $N_p$, and $\lambda^i_p$ are the eigenvalues returned by the POD.
}
{Figures \ref{reduced velocity non parametric} and \ref{reduced displacement} show two reproductive reduced order solutions: the fluid velocity and the solid displacement, respectively; as we can see from Figure \ref{reduced displacement}, the reduced order model shows a good capability also in reproducing very small deformations: these results were obtained using $N_z=15$, $N_p=10$ and $N_d=10$ modes, respectively. Figure \ref{error behavior} shows that, with $N_z=15$, $N_p =10$ and $N_d=10$ basis functions for each component of the solution, we have a good relative approximation error behavior over time: as we can see from the figure, at the final timesteps of the simulation the error increases, and we think this is due to some error accumulation phenomenon. The error has been computed as the $L^2$ error for the fluid pressure, and as the $H^1$ error for the fluid velocity and the solid displacement:
{
\begin{equation*}
\text{err}_p(t^{i+1}) = \frac{\lvert\lvert p_{f, h}^{i+1} - p_{f, {N_p}}^{i+1}\rvert\rvert_{L^2(\Omega_f)}}{\lvert\lvert p_{f, h}^{i+1}\rvert\rvert_{L^2(\Omega_f)}},
\end{equation*}
\begin{equation*}
\text{err}_d(t^{i+1}) = \frac{\lvert\lvert \bm{d}_{s, h}^{i+1} - \bm{d}_{s, {N_d}}^{i+1}\rvert\rvert_{H^1(\Omega_s)}}{\lvert\lvert \bm{d}_{s, h}^{i+1}\rvert\rvert_{H^1(\Omega_s)}}.
\end{equation*}
\begin{equation*}
\text{err}_u(t^{i+1}) = \frac{\lvert\lvert \bm{u}_{f, h}^{i+1} - \bm{u}_{f, {N_z}}^{i+1}\rvert\rvert_{H^1(\Omega_f)}}{\lvert\lvert \bm{u}_{f, h}^{i+1}\rvert\rvert_{H^1(\Omega_f)}}.
\end{equation*}
}
We were interested in seeing how the average approximation error and the average number of internal iterations changes, by changing the number of reduced basis $N_p$ used for the fluid pressure in the reduced order model: results are reported in Table \ref{table average relative error non parametric} and Table \ref{table average internal iterations non parametric} respectively. {As we can see from Table \ref{table average relative error non parametric}, the average approximation error decreases up to when we use $N_p = 25$ modes for the pressure, then we observe an increment in the approximation error: we read this result as the fact that with $40$ modes for the fluid pressure, we are just adding noise to the online system. It is also interesting to see from Table \ref{table average internal iterations non parametric} that the average number of internal iterations required from the algorithm, in order to reach a coupling tolerance of $\varepsilon = 10^{-6}$ is relatively higher for the reduced order model, when compared to the full order one: this is due to the reduction of the two problems, but in any case we can see that this number stabilizes around $26-27$ subiterations.}
{Finally, Figure \ref{internal iterations behavior} depicts the behavior with respect to time of the number of internal iterations: a comparison is drawn between the full order model and the reduced order one, where we used $N_p=10$ reduced basis functions. We can see that the number of internal iterations, both for the offline and for the online part stabilizes towards the end of the simulation.} 
We would like to make the following remark: all these results are computed by varying the number of modes used for the approximation of the fluid pressure, while keeping fixed both $N_{\bm{z}}$ and $N_{\bm{d}}$. {The motivation behind our choice is the fact that we want to see how the number of modes directly impacts the performance of the method, and, more precisely, of the implicit step, where the coupling between the two physics is imposed, by coupling the pressure Poisson problem with the solid problem.}}
The authors are aware that these results are by no means exhaustive, and this is a further testimony to the capability that such a partitioned procedure offers: many more tests are possible, where for example $N_{\bm{d}_s}$ is varied, and $N_{\bm{z}}$ is kept fixed, or both can vary.

\begin{figure}
\centering
\includegraphics[scale=0.3]{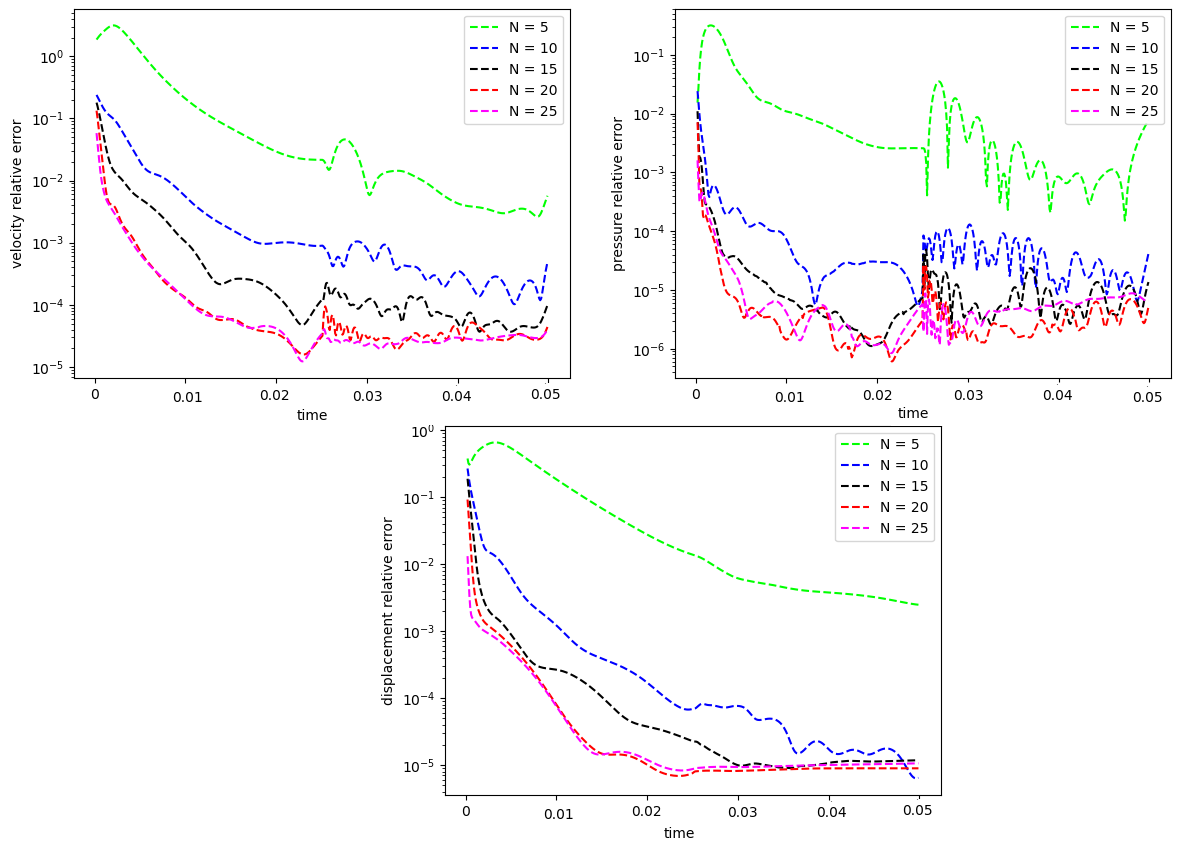}
\caption{{Relative error behavior in time, by varying the number $N$ of modes used in the online phase. Top left: fluid velocity error. Top right: fluid pressure error. Bottom center: solid displacement error.}}
\label{error in time same N}
\end{figure}

\begin{figure}
\centering
\includegraphics[scale=0.3]{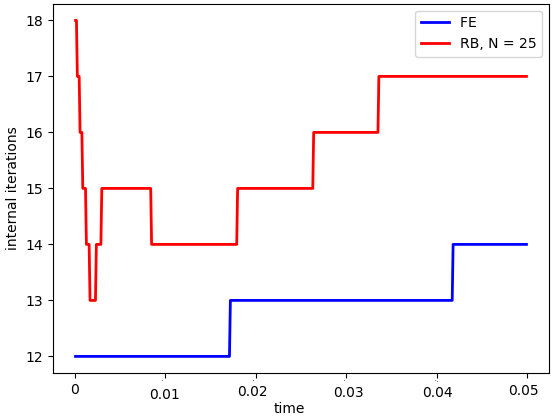}
\caption{{Internal iterations: comparison between full order model (blue) and reduced order model (red). The reduced solutions have been obtained with: $N = 25$ modes.}}
\label{iterations in time same N}
\end{figure}

\begin{table}
\centering
\caption{{Average number of subiterations: comparison between full order model, and reduced order model with basis refinement for all the components $\bm{u}_f$, $p_f$ and $\bm{d}_s$.}}
\label{table average internal iterations same N}
\begin{tabular}{|l | c| c|}
\hline
$N$ & full order model & reduced order model\\
\hline
$5$ & $13$ & $17$\\ 
$10$ & $13$ & $25$\\
$15$ & $13$ & $25$\\
$20$ & $13$ & $24$\\
$25$ & $13$ & $16$\\
\hline
\end{tabular}
\end{table}

{For this reason, we now present some additional results, that have been instead obtained by using the same number of modes for all the components of the FSI solution, namely $N_z=N_p=N_d=N$. Figure \ref{error in time same N} shows the relative approximation error in time, for the three components of the FSI solution: as we can see, by increasing the number of modes, we get a better approximation error. In particular, from Figure \ref{error in time same N}, we observe an oscillating behavior of the pressure relative error towards the last timesteps of the simulation: this behavior can be seen also in Figure \ref{error behavior}, and can be interpreted as the result of an accumulation phenomenon, where the error accumulates and starts to oscillate. Another important observation is that the online pressure, with this partitioned approach, has been obtained without the supremizer enrichment technique \cite{Supremizer}; as already remarked in \cite{Ballarin2017}, this may lead to non optimal error convergence. Motivated by this, we think it can be a very interesting point for a future work to see if this has some implications in the error oscillation that we observe in Figure \ref{error in time same N}. Figure \ref{iterations in time same N} shows the number of internal iterations for each timestep of the simulation: the online computations are performed using $N = 25$ modes for all the components, and are compared to the high fidelity computations. Also in this case the average number is higher for the online simulation, due to model order reduction, and also in this case the number stabilizes around $17$ iterations towards the end of the simulation. Finally, in Table \ref{table average internal iterations same N} we report the average number of internal iterations, for an increasing number $N$ of modes used in the online phase. As we can see, by increasing $N$ from $10$ to $25$ there is almost no improvement in the number of internal iterations: it stabilizes around $25$. This number then drops to $16$ (which is very close to the FOM results) for $N=25$: we did not increase $N$ further, because $N_d=25$ is the total number of modes retained by the POD on the solid displacement, before hitting a very small magnitude (less than $10^{-9}$) for the corresponding eigenvalues.}
\section{Shape parametrization of the leaflets}\label{geometrical parametrization}
In this section we are going to address a slightly different situation, the difference being now the presence of a geometrical parameter $\mu_g$, that represents the length of the leaflets; we also admit the possibility of a further physical parameter $\mu_p$, so that,  to summarize, we consider a parameter $\mu\in\mathcal{P}\subset\mathbb{R}^d$, where $d=1$ if just a geometrical parametrization is considered (and thus $\mu=\mu_g$), or $d=2$ (and thus $\mu=(\mu_g, \mu_p)$).
\subsection{FSI in the presence of shape parametrization}\label{ALE geometrical formalism}
Let us denote by $\Omega(t;\mu_g):=\Omega_f(t;\mu_g)\cup\Omega_s(t;\mu_g)$ the current physical domain: we now have a time dependence and a parameter dependence. We introduce the time--independent \emph{intermediate configuration} $\tilde{\Omega}(\mu_g):=\tilde{\Omega}_f(\mu_g)\cup\tilde{\Omega}_s(\mu_g)$, where we are considering the reference configuration of both physics, still taking into account the parameter dependence. Finally, we have the time--independent, parameter--independent \emph{reference configuration} $\hat{\Omega}:=\hat{\Omega}_f\cup\hat{\Omega}_s$.
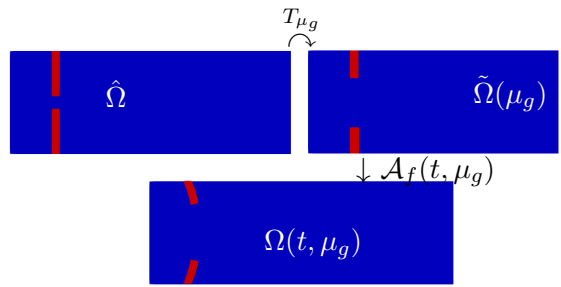
\begin{figure}
\centering
\begin{tikzpicture}
\node[anchor=south west,inner sep=0] (image) at (0,0) {\includegraphics[scale=0.2]{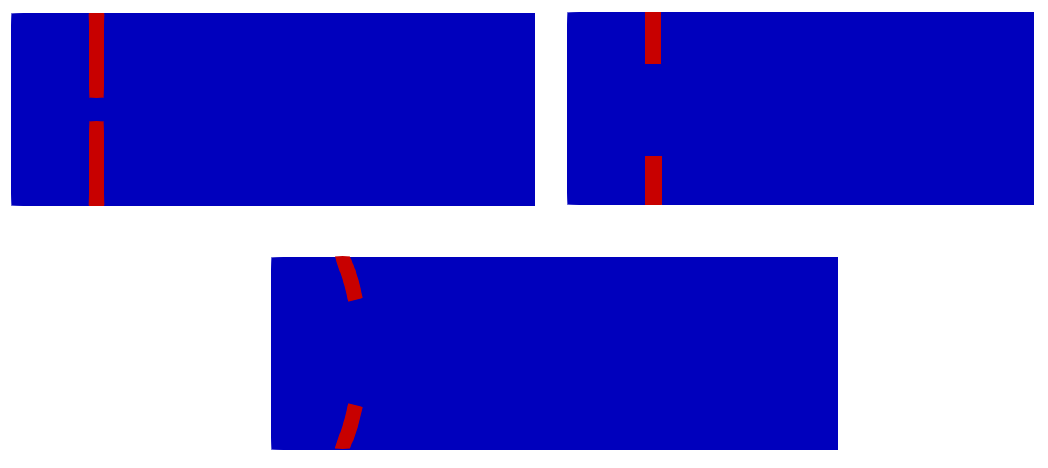}};
\begin{scope}[x={(image.south east)},y={(image.north west)}]
	\draw[white,thick] (0.9,0.8) node {$\tilde{\Omega}(\mu_g)$};
	\draw[white,thick] (0.2,0.8) node {$\hat{\Omega}$};
	\draw[black, thick] (0.53, 1.07) node {$\overset{T_{\mu_g}}{\curvearrowright}$};
	\draw[white, thick] (0.55, 0.2) node {$\Omega(t, \mu_g)$};
	\draw[black, thick] (0.75, 0.49) node {$\downarrow$ $\mathcal{A}_f(t, \mu_g)$};
\end{scope}
\end{tikzpicture}
\caption{Domains: reference configuration $\hat{\Omega}$ (top left), parametrized reference configuration $\tilde{\Omega}(\mu_g)$ (top right), and original configuration $\Omega(t;\mu_g)$ (bottom).}
\label{domain parametrization}
\end{figure}

We call $T$ the shape parametrization map; for every $\mu_g$ we have a map $T_{\mu_g}$ defined as follows:
\begin{equation*}
\begin{split}
T_{\mu_g}\colon\hat{\Omega}&\mapsto\tilde{\Omega}(\mu_g)\\
\hat{x}&\mapsto\tilde{x}=T_{\mu_g}(\hat{x}).
\end{split}
\end{equation*}
We then have the ALE map $\mathcal{A}_f(t;\mu)$, already introduced in Section \ref{ALE formalism}, which is now a map from the current parametrized fluid configuration $\hat{\Omega}_f(t;\mu_g)$ and the intermediate fluid configuration $\tilde{\Omega}_f(\mu_g)$:
\begin{equation*}
\begin{split}
\mathcal{A}_f(t;\mu)\colon\tilde{\Omega}_f(\mu_g)&\mapsto\hat{\Omega}_f(t;\mu_g)\\
\tilde{x}&\mapsto\hat{x} = \tilde{x}+\tilde{\bm{d}}_f(\tilde{x}; t, \mu),
\end{split}
\end{equation*}
where $\tilde{\bm{d}}_f$ is the mesh displacement already defined in Section \ref{ALE formalism}.
\\ Let us define the gradients and the determinants of the deformation maps:
\begin{gather}
\bm{G}(\hat{x}; \mu_g) = \hat{\nabla} T_{\mu_g}(\hat{x}), \qquad K(\hat{x}; \mu_g) = \text{det}\bm{G}(\hat{x}; \mu_g),\nonumber\\
\tilde{\bm{F}}(\tilde{x};\mu) = \tilde{\text{Id}} + \tilde{\nabla}\tilde{\bm{d}}_f(\mu), \qquad \tilde{J}(\tilde{x}; \mu) =\text{det}\tilde{\bm{F}}.
\end{gather}
We can pull-back the gradient $\tilde{\bm{F}}(\tilde{x}; \mu)$ to the reference domain $\hat{\Omega}_f$, and we obtain $\bm{F}(\hat{x}; \mu) = \text{Id} + \hat{\nabla}\hat{\bm{d}}_f(\mu)\bm{G}^{-1}(\hat{x}, \mu_g)$. With this notation, we can conclude that the gradient of the deformation map from the reference configuration to the current configuration is given by $\bm{F}(\hat{x}, \mu)\bm{G}(\hat{x}, \mu_g)$; let us denote by $\tilde{\bm{F}}_{\mu}$, $\bm{F}_{\mu}$ and $\bm{G}_{\mu_g}$ the gradients $\tilde{\bm{F}}(\tilde{x}, \mu)$, $\bm{F}(\hat{x}, \mu)$ and $\bm{G}(\hat{x}, \mu_g)$ respectively, and by $J_{\mu}$ and $K_{\mu_g}$ the determinants of $\bm{F}_{\mu}$ and $\bm{G}_{\mu_g}$. We are now ready to state the strong form of the problem of interest.
\subsection{Strong formulation}\label{parametric problem formulation}
The strong form of the parametrized FSI problem reads as follows: for every $t\in[0, T]$ and for every $\mu\in\mathcal{P}$, find the fluid velocity $\bm{u}_f(t;\mu)\colon \Omega_f(t;\mu_g)\mapsto\mathbb{R}^2$, the fluid pressure $p_f(t;\mu)\colon\Omega_f(t;\mu_g)\mapsto\mathbb{R}$, the mesh displacement $\tilde{\bm{d}}_f(t;\mu)\colon\tilde{\Omega}_f(\mu_g)\mapsto\mathbb{R}^2$ and the solid displacement $\tilde{\bm{d}}_s(t;\mu)\colon\tilde{\Omega}_s(\mu_g)\mapsto\mathbb{R}^2$ such that:
\begin{equation*}
\begin{cases}
-\tilde{\Delta}\tilde{\bm{d}}_f = 0 \quad\text{in }\tilde{\Omega}_f(\mu_g)\times[0, T],\\
\tilde{\bm{d}}_f = \tilde{\bm{d}}_s \quad\text{on }\tilde{\Gamma}_{FSI}(\mu_g)\times[0, T],
\end{cases}
\end{equation*}
and
\begin{equation*}
\begin{cases}
\rho_f\partial_t\bm{u}_f\lvert_{\tilde{x}} + \rho_f(\bm{u}_f - \partial_t\bm{d}_f\lvert_{\tilde{x}})\cdot \nabla\bm{u}_f - \text{div}\sigma_f(\bm{u}_f, p_f) =0 \quad\text{in }\Omega_f(t;\mu_g)\times[0, T],\\
\text{div}\bm{u}_f = 0 \quad\text{in }\Omega_f(t;\mu_g)\times[0, T],\\
\rho_s\partial_{tt}\tilde{\bm{d}}_s - \tilde{\text{div}}(\tilde{\bm{P}}(\tilde{\bm{d}}_s)) =0 \quad\text{in }\tilde{\Omega}_s(\mu_g)\times[0, T].
\end{cases}
\end{equation*}
Here we notice that, again, the fluid problem is formulated in the current parametrized configuration $\Omega_f(t;\mu_g)$, whereas the solid problem is formulated in the parametrized intermediate configuration $\tilde{\Omega}_s(\mu_g)$. The quantity $\partial_t\bm{u}_f\lvert_{\tilde{x}}$ represents the ALE time derivative: $\partial_t\bm{u}_f(\bm{x}, t;\mu_g)\lvert_{\tilde{x}}= \partial_t\tilde{\bm{u}}_f(\tilde{\bm{x}}, t;\mu_g)$. Again, $\sigma_f$ is the fluid Cauchy stress tensor, and $\tilde{\bm{P}}$ is the first Piola--Kirchoff stress tensor: their definition has been given in Section \ref{problem formulation}.
The previous system is completed by some suitable initial conditions, {by the same boundary conditions prescribed in \eqref{boundary conditions}} and by the following coupling conditions:
\begin{equation*}
\begin{cases}
\bm{u}_f = \frac{d}{dt}\bm{d}_s \text{ on }\Gamma_{FSI}(t;\mu_g),\\
\tilde{\bm{d}}_f = \tilde{\bm{d}}_s \text{ on }\tilde{\Gamma}_{FSI}(\mu_g),\\ 
\tilde{J}_{\mu_g}\tilde{\sigma_f}(\tilde{\bm{u}}_f, \tilde{p}_f)\tilde{\bm{F}}^{-T}_{\mu_g}\tilde{\bm{n}}_f = -\tilde{\bm{P}}(\tilde{\bm{d}}_s)\tilde{\bm{n}}_s \text{ on }\tilde{\Gamma}_{FSI}(\mu_g),
\end{cases}
\end{equation*}
being $\tilde{\sigma}_f$ the Cauchy stress tensor in the parametrized intermediate fluid domain $\tilde{\Omega}_f(\mu_g)$:
\begin{equation*}
\tilde{\sigma}_f(\tilde{\bm{u}}_f, \tilde{p}_f) = \mu_f(\tilde{\nabla}\tilde{\bm{u}}_f\tilde{\bm{F}}^{-1}_{\mu} + \tilde{\bm{F}}^{-T}_{\mu}\tilde{\nabla}^T\tilde{\bm{u}}_f).
\end{equation*}
Thanks to the introduction of the pull--back maps, we can reformulate our problem in the reference configuration $\hat{\Omega}$: for every $t\in [0, T]$ and for every $\mu\in\mathcal{P}$, find the fluid velocity $\hat{\bm{u}}_f(t,\mu):\hat{\Omega}_f\mapsto\mathbb{R}^2$, the fluid pressure $\hat{p}_f(t,\mu):\hat{\Omega}_f\mapsto\mathbb{R}$, the fluid displacement $\hat{\bm{d}}_f(t,\mu)\colon\hat{\Omega}_f\mapsto\mathbb{R}^2$ and the solid deformation $\hat{\bm{d}}_s(t, \mu):\hat{\Omega}_s\mapsto\mathbb{R}^2$ such that:
\begin{equation}
\begin{cases}
\begin{split}
&\rho_fJ_{\mu}K_{\mu_g}(\partial_t\hat{\bm{u}}_f + \hat{\nabla}\hat{\bm{u}}_f\bm{G}^{-1}_{\mu_g}\bm{F}^{-1}_{\mu}(\hat{\bm{u}}_f-\partial_t\hat{\bm{d}}_f))- \hat{\text{div}}(J_{\mu}K_{\mu_g}\hat{\sigma}_f(\hat{\bm{u}}_f, \hat{p}_f)\bm{F}^{-T}_{\mu}\bm{G}^{-T}_{\mu_g}) = \\
&=0 \quad\text{in }\hat{\Omega}_f\times(0, T],
\end{split}
\\
\hat{\text{div}}(J_{\mu}K_{\mu_g}\bm{G}^{-1}_{\mu_g}\bm{F}^{-1}_{\mu}\hat{\bm{u}}_f) = 0 \quad \text{in }\hat{\Omega}_f\times(0, T],\\
-\hat{\text{div}}(K_{\mu_g}\hat{\nabla}\hat{\bm{d}}_f\bm{G}^{-1}_{\mu_g}\bm{G}^{-T}_{\mu_g}) = 0 \quad \text{in }\hat{\Omega}_f\times(0, T],\\
\rho_sK_{\mu_g}\partial_{tt}\hat{\bm{d}_s} -\hat{\text{div}}(K_{\mu_g}\hat{\bm{P}}(\hat{\bm{d}_s})\bm{G}^{-T}_{\mu_g}) = 0\quad \text{in }\hat{\Omega}_s\times(0, T],
\end{cases}
\end{equation}
where:
\begin{gather}
\hat{\sigma}_f(\hat{\bm{u}}_f, \hat{p}_f) = \mu_f(\hat{\nabla}\hat{\bm{u}}_f\bm{G}^{-1}_{\mu_g}\bm{F}^{-1}_{\mu} + \bm{F}^{-T}_{\mu}\bm{G}^{-T}_{\mu_g}\hat{\nabla}^T\hat{\bm{u}}_f),\nonumber\\
\hat{\bm{P}}(\hat{\bm{d}}_s) = \lambda_s\text{tr}\varepsilon_s(\hat{\bm{d}}_s)\bm{I} + 2\mu_s\varepsilon_s(\hat{\bm{d}}_s),\nonumber\\
\varepsilon_s(\hat{\bm{d}}s) = \frac{1}{2}(\hat{\nabla}\hat{\bm{d}}_s\bm{G}^{-1}_{\mu_g} + \bm{G}^{-T}_{\mu_g}\hat{\nabla}^T\hat{\bm{d}}_s).
\end{gather}
We have the coupling conditions
\begin{equation}
\label{coupling conditions}
\begin{cases}
\hat{\bm{d}}_f = \hat{\bm{d}}_s \quad\text{on }\hat{\Gamma}_{FSI}\\ 
\hat{\bm{u}}_f = \partial_t\hat{\bm{d}}_s\quad \text{on }\hat{\Gamma}_{FSI},\\
J_{\mu_g}K_{\mu_g}\hat{\sigma_f}(\hat{\bm{u}}_f, \hat{p}_f)\bm{F}^{-T}_{\mu_g}\bm{G}^{-T}_{\mu_g}\hat{\bm{n}}_f = -K_{\mu_g}\hat{\bm{P}}(\hat{\bm{d}}_s)\bm{G}^{-T}_{\mu_g}\hat{\bm{n}}_s \quad\text{on }\hat{\Gamma}_{FSI},
\end{cases}
\end{equation}
and the following boundary conditions:
\begin{equation*}
\begin{cases}
\hat{\sigma}_f(\hat{\bm{u}}_f, \hat{p}_f)\hat{\bm{n}} = -p_{in}(t)\hat{\bm{n}}\text{ on $\hat{\Gamma}_{in}$},\\
\hat{\sigma}_f(\hat{\bm{u}}_f, \hat{p}_f)\hat{\bm{n}} = -p_{out}(t)\hat{\bm{n}}\text{ on $\hat{\Gamma}_{out}$},\\
\hat{\bm{d}}_s = 0\text{ on $\hat{\Gamma}_s^D$},\\
\hat{\bm{u}}_f = 0\text{ on }\hat{\Gamma}_{top}\cup \hat{\Gamma}_{bott}.
\end{cases}
\end{equation*}
Again, $\hat{\bm{n}}$ represents the normal vector to the relative part of the boundary of the domain.
\begin{remark}
In this section we stressed the difference between entities on the current configuration, the parametrized intermediate configuration and the reference configuration, by using the superscripts $\hat{}$ and $\tilde{}$. However, since from now on everything will be cast in the reference configuration, and in order to make the notation as light as possible, we will drop all the superscripts.
\end{remark}
\subsection{Offline computational phase}\label{offline phase parametrized}
Hereafter we present the offline phase of the partitioned procedure in the presence of a parameter. We employ again a Chorin-Temam projection scheme for the Navier--Stokes equation; we define a time-stepping procedure by sampling the time interval $[0, T]$ with an equispaced sampling $\{t_0, \hdots, t_{N_T}\}$, where $t_i = i\Delta T$, for $i=0, \hdots, N_T$ and $N_T=\frac{T}{\Delta T}$. We discretize the time derivative of a function $f$ with a first backward difference: $D_tf^{i+1} = \frac{f^{i+1} - f^i}{\Delta T}$, and $D_{tt}f^{i+1}= D_t(D_tf^{i+1})$, where $f^{i+1}=f(t^{i+1})$.  We then discretize the parameter space $\mathcal{P}\subset\mathbb{R}^d$, $d=1,2$ with an equispaced sampling, and we obtain $\mathcal{P}_{train}=[\mu_g^1,\hdots,\mu_g^{N_g}]\times[\mu_s^1,\hdots, \mu_s^{N_s}]$ (in the case of physical and geometrical parametrization) or $\mathcal{P}_{train}=[\mu_g^1,\hdots,\mu_g^{N_g}]$ (in the case of just geometrical parametrization). We define here $N_{train}$ to be the cardinality of the training set $\mathcal{P}_{train}$.\\
In the following, we use the same function spaces that we have introduced in Section \ref{space discretization section}:
\begin{gather} 
V(\Omega_f) := [H^1(\Omega_f)]^2,\nonumber\\
E^f(\Omega_f) := [H^1(\Omega_f)]^2,\nonumber\\
Q(\Omega_f):= L^2(\Omega_f),\nonumber\\
E^s(\Omega_s)= [H^1(\Omega_s)]^2,\nonumber
\end{gather}
endowed with the $H^1$ norm ($V(\Omega_f)$), the $H^1$ seminorm ($E^f(\Omega_f)$ and $E^s(\Omega_s)$) and the $L^2$ norm respectively. We remark that in the previous definitions, the domains $\Omega_f$ and $\Omega_s$
are the \emph{reference configurations} (both parameter-- and time-- independent).
Again we discretize in space the FSI problem, using second order Lagrange Finite Elements for fluid velocity resulting in the discrete space $V_h\subset V$, while the fluid pressure, the fluid displacement and the solid displacement are discretized with first order Lagrange Finite Elements, resulting in the discrete spaces $Q_h\subset Q$, $E^f_h\subset E^f$ and $E_h^s\subset E^s$; we make again use of a lifting function for the fluid pressure, thus we introduce also the discrete space $Q_h^0$, which is defined exactly as in Section \ref{space discretization section}. 
\begin{remark}
In this case the lifting function $\ell(t)$ does not depend on the parameter $\mu_g$, as we can deduce from the fact that the quantity $p_{in}(t)$ is parameter--independend. Therefore we can compute the lifting function, during the offline phase, \emph{once and for all} for every timestep $t^i$. {We would also like to refer the interested reader to the work presented in \cite{GUNZBURGER20071030} in the case $p_{in}$ is parameter dependent: indeed, in \cite{GUNZBURGER20071030} the authors present a detailed description of the work that has to be done in order to implement reduced order models in the presence of multiple parameters in the boundary data.}
\end{remark}
The space discretized version of the partitioned procedure now reads as follows: for $i=0, \hdots, N_T$, for $\mu=(\mu_g, \mu_p)\in\mathcal{P}_{train}$:
\subsubsection*{\underline{Extrapolation of the mesh displacement}:} 
find $\bm{d}_{f, h}^{i+1}\in E_h^f$ such that $\forall \bm{e}_{f,h} \in E^f_h$:
\begin{equation}
\label{mesh displacement}
\begin{cases}
\int_{\Omega_f}K_{\mu_g}\nabla\bm{d}_{f, h}^{i+1}\bm{G}^{-1}_{\mu_g}\cdot\nabla\bm{e}_{f, h}\bm{G}^{-1}_{\mu_g}\,dx =0 \\
\bm{d}_{f, h}^{i+1} = \bm{d}_{s, h}^i \quad\text{on $\Gamma_{FSI}$}. 
\end{cases}
\end{equation}
\subsubsection*{\underline{Fluid explicit step}:}
find $\bm{u}_{f, h}^{i+1}\in V_h$ such that $\forall \bm{v}_h\in V_h$:
\begin{equation}
\label{fluid explicit step}
\begin{cases}
\begin{split}
&\rho_f\int_{\Omega_f} J_{\mu}K_{\mu_g}\Bigl(\frac{\bm{u}_{f, h}^{i+1}-\bm{u}_{f, h}^i}{\Delta T}\Bigr)\cdot\bm{v}_h\,dx + \\
&+ \rho_f\int_{\Omega_f}J_{\mu}K_{\mu_g}[\nabla\bm{u}_{f, h}^{i+1}\bm{G}^{-1}_{\mu_g}\bm{F}^{-1}_{\mu}](\bm{u}_{f, h}^{i+1} -D_t\bm{d}_{f, h}^{i+1})\cdot\bm{v}_h\,dx\\
&+ \mu_f\int_{\Omega_f}J_{\mu}K_{\mu_g}\varepsilon(\bm{u}_{f, h}^{i+1})\bm{F}^{-T}_{\mu}\bm{G}^{-T}_{\mu_g}:\nabla\bm{v}_h\,dx + \int_{\Omega_f}J_{\mu}K_{\mu_g}\bm{F}^{-T}_{\mu}\bm{G}^{-T}_{\mu_g}\nabla p_{f, h}^i\cdot\bm{v}_h\,dx= 0
\end{split}
\\ \bm{u}_{f, h}^{i+1} = D_t\bm{d}_{f, h}^{i+1} \quad \text{on $\Gamma_{FSI}$.}
\end{cases}
\end{equation}
\subsubsection*{\underline{Implicit step}:} for any $j=0, \dots$ until convergence:
\begin{enumerate}
\item \textbf{fluid projection substep (pressure Poisson formulation):} find $p_{f, h}^{0, i+1, j+1}\in Q_h$ such that $\forall q_h\in Q_h^0$:
\begin{equation*}
\begin{split}
&\alpha_{ROB}\int_{\Gamma_{FSI}}p_{f, h}^{0, i+1, j+1}q_h\,ds +  \int_{\Omega_f}J_{\mu}K_{\mu_g}\bm{F}^{-T}_{\mu}\bm{G}^{-T}_{\mu_g}\nabla p_{f, h}^{0, i+1, j+1}\cdot \bm{F}^{-T}_{\mu}\bm{G}^{-T}_{\mu_g}\nabla q_h\,dx =\\
&-\frac{\rho_f}{\Delta T}\int_{\Omega_f}\text{div}(J_{\mu}K_{\mu_g}\bm{F}^{-1}_{\mu}\bm{G}^{-1}_{\mu_g}\bm{u}_{f, h}^{i+1})q_h\,dx + \alpha_{ROB}\int_{\Gamma_{FSI}}p_{f, h}^{i+1, j}q_h\,ds  - \\
&-\alpha_{ROB}\int_{\Omega_f}\ell^{i+1}\cdot q_h\,dx - \rho_f\int_{\Gamma_{FSI}}(D_{tt}\bm{d}_{s, h}^{i+1, j})\cdot J_{\mu}K_{\mu_g}\bm{F}^{-T}_{\mu}\bm{G}^{-T}_{\mu_g}\bm{n}_fq_h\,ds -\\
& -\int_{\Omega_f}J_{\mu}K_{\mu_g}\bm{F}^{-T}_{\mu}\bm{G}^{-T}_{\mu_g}\nabla \ell^{i+1}\cdot\bm{F}^{-T}_{\mu}\bm{G}^{-T}_{\mu_g}\nabla q_h\,dx
\end{split}
\end{equation*}
subject to the boundary conditions \eqref{pressure BC}. We then retrieve the original fluid pressure $p_{f, h}^{i+1, j+1} = p_{f,h}^{0, i+1, j+1}+\ell^{i+1}$.
\item \textbf{structure projection substep:} find $\bm{d}_{s, h}^{i+1, j+1}\in E^s_h$ such that $\forall \bm{e}_{s, h}\in E^s_h$:
\begin{equation*}
\begin{split}
&\rho_s\int_{\Omega_s}K_{\mu_g}D_{tt}\bm{d}_{s, h}^{i+1, j+1}\cdot\bm{e}_{s, h}\,dx  + \int_{\Omega_s}K_{\mu_g}\bm{P}(\bm{d}_{s, h}^{i+1, j+1})\bm{G}^{-T}_{\mu_g}:\nabla\bm{e}_{s, h}\,dx= \\
& = -\int_{\Gamma_{FSI}}J_{\mu_g}K_{\mu_g}\sigma_f(\bm{u}_{f, h}^{i+1, j+1}, p_{f, h}^{i+1, j+1})\bm{F}^{-T}_{\mu_g}\bm{G}^{-T}_{\mu_g}\bm{n}_f\cdot\bm{e}_{s, h}\,dx
\end{split}
\end{equation*}
subject to the boundary condition $\bm{d}_{s,h}^{i+1, j+1} = 0$ on $\Gamma_{D}^s$.
\end{enumerate}
In the fluid projection step, in order to enhance the stability of the method we have employed again a Robin boundary condition, which in the case of shape parametrization reads as follows:
\begin{equation*}
\begin{split}
\alpha_{ROB}p^{i+1} + \bm{F}^{-T}_{\mu}\bm{G}^{-T}_{\mu_g}\nabla p^{i+1}&\cdot J_{\mu}K_{\mu_g}\bm{F}^{-T}_{\mu_g}\bm{G}^{-T}_{\mu_g}\bm{n}_f = \\
&= \alpha_{ROB} p^{i+1,\star} -\rho_fD_{tt}\bm{d}_s^{i+1,\star}\cdot J_{\mu}K_{\mu_g}\bm{F}^{-T}_{\mu}\bm{G}^{-T}_{\mu_g}\bm{n}_f.
\end{split}
\end{equation*}
\subsubsection{POD--Greedy}
For this test case we decided to adopt a POD strategy which is slightly different with respect to the standard POD that we presented for the other two test cases. Indeed, the idea here is to first perform a standard POD in time on each set of snapshots computed for each value of the parameter ${\mu}$ in the training set $\mathcal{P}_{train}$. Then, we take all the modes computed with the standard POD, weighted with the corresponding eigenvalue, and perform a final outer run of POD. {The idea that we implement here is inspired by the POD--Greedy strategy presented in \cite{Nguyen}, and the motivation behind our choice is given by the fact that now the parameter space has a higher dimension (two parameters are being considered): for this reason, performing a standard POD would be computationally unfeasible. Indeed we work with a huge number of collected snapshots, which result in a correlation matrix of enormous dimension: solving the eigenvalue problem for this matrix quickly saturates the RAM of a computer; therefore performing a ``naive'' POD is not a good idea in this case.}
Let us now present briefly the procedure: we consider all the parameters $\mu_i\in\mathcal{P}_{train}$ in the training set: here the index $i$ has to be considered as a single index, in the case of geometrical parametrization only ($\mu=\mu_g$), or as a multiindex $i=(i_g, i_p)$ in the case of geometrical and physical parametrization ($\mu=(\mu_g, \mu_p)$).
We start by constructing, for each parameter ${\mu}_i\in \mathcal{P}_{train}$, the snapshots matrices $\bm{\mathcal{S}}_z({\mu}_i)$ for the fluid change of variable $\bm{z}_f$, $\bm{\mathcal{S}}_p(\mu_i)$ for the fluid pressure $p_f^0$ and $\bm{\mathcal{S}}_{d_s}({\mu}_i)$ for the solid displacement $\bm{d}_s$:
\begin{equation*}
\begin{split}
\bm{\mathcal{S}}_z({\mu}_{i})=\{\bm{z}_{f, h}(t_0, {\mu}_{i}), \hdots, \bm{z}_{f, h}(t_{N_T}, {\mu}_{i})\} \in \mathbb{R}^{\mathcal{N}_u^h\times N_T},\\
\bm{\mathcal{S}}_p({\mu}_{i})=\{p_{f, h}^0(t_0, {\mu}_{i}),\hdots, p_{f, h}^0(t_{N_T}, {\mu}_{i})\} \in \mathbb{R}^{\mathcal{N}_p^h\times N_T},\\
\bm{\mathcal{S}}_{d_s}({\mu}_{i})=\{\bm{d}_{s,h}(t_0, {\mu}_{i}), \hdots, \bm{d}_{s, h}(t_{N_T}, {\mu}_{i})\} \in \mathbb{R}^{\mathcal{N}_{d_s}^h\times N_T}.
\end{split}
\end{equation*}
We then perform a standard POD on each snapshots matrix and we extract the basis functions $\{\bm{\Phi}_{z_f}^k({\mu}_{i})\}_{k=1}^{N^z_{i}}$, $\{\Phi_p^k({\mu}_{i})\}_{k=1}^{N^p_{i}}$ and $\{\bm{\Phi}_{d_s}^k({\mu}_{i})\}_{k=1}^{N^d_{i}}$. Let us also call $\{\lambda_n^z\}_{n=1}^{N^z_{i}}$, $\{\lambda_n^p\}_{n=1}^{N^p_{i}}$, $\{\lambda_n^d\}_{n=1}^{N^d_{i}}$ the eigenvalues, ordered by decreasing order of magnitude, returned by the POD on each snapshot matrix $\bm{\mathcal{S}}_z({\mu}_{i})$, $\bm{\mathcal{S}}_p({\mu}_{i})$, $\bm{\mathcal{S}}_{d_s}({\mu}_{i})$.
Afterwards, we perform a second run of POD in the following way: we start by building the snapshots matrices always for the components $\bm{z}_f$, $p_f$ and $\bm{d}_s$, weighting each snapshot with the corresponding eigenvalue given by the standard POD:
\begin{equation*}
\begin{split}
\bm{\mathcal{S}}_z=\{\sqrt{\lambda_1^z}\bm{\Phi}_{z_f}^1({\mu}_{1}), \hdots, {\sqrt{\lambda^z_{N^z_{1}}}}\bm{\Phi}_{z_f}^{N^z_{1}}({\mu}_{1}), \hdots, {\sqrt{\lambda^z_{N^z_{N_{train}}}}}\bm{\Phi}_{z_f}^{N^z_{N_{train}}}({\mu}_{N_{train}})\},\\
\bm{\mathcal{S}}_p=\{{\sqrt{\lambda_1^p}}\Phi_p^1({\mu}_{1}), \hdots, {\sqrt{\lambda^p_{N^p_{1}}}}\Phi_p^{N^p_{1}}({\mu}_{1}), \hdots, {\sqrt{\lambda^p_{N^p_{N_{train}}}}}\Phi_p^{N^p_{N_{train}}}({\mu}_{N_{train}})\},\\
\bm{\mathcal{S}}_{d_s}=\{{\sqrt{\lambda_1^d}}\bm{\Phi}_{d_s}^1({\mu}_{1}), \hdots, {\sqrt{\lambda^d_{N^d_{1}}}}\bm{\Phi}_{d_s}^{N^d_{1}}({\mu}_{1}), \hdots, {\sqrt{\lambda^d_{N^d_{N_{train}}}}}\bm{\Phi}_{d_s}^{N^d_{N_{train}}}({\mu}_{N_{train}})\}.
\end{split}
\end{equation*}
{The weighting of the POD modes obtained from the first POD run is motivated by the fact that, in this way, the second POD will be correctly weighted to accommodate modes from different parameter values: for a detailed discussion, we refer the interested reader to \cite{Nguyen}.}
We the perform a second POD on the previous snapshots matrices, and we finally obtain a set of basis functions $\{\bm{\Phi}_{z_f}^k\}_{k=1}^{N^z}$, $\{\Phi_p^k\}_{k=1}^{N^p}$ and $\{\bm{\Phi}_{d_s}^k\}_{k=1}^{N^d}$. Then, to obtain a set of basis functions for the mesh displacement $\bm{d}_f$, we choose to employ again an harmonic extension of the solid displacement basis functions $\bm{\Phi}_{d_s}^k$ on the entire fluid domain $\Omega_f$, as we have done for the previous two test cases; in this way we obtain a set $\{\bm{\Phi}_{d_f}^k\}_{k=1}^{N^d}$ of reduced basis also for the mesh displacement.
\subsection{Online phase}\label{shape and physical parametrization: online phase}
We start by introducing the online solutions $\bm{z}_{f, N}^{i+1}(\mu)$, $p_{f, N}^{i+1}(\mu)$ and $\bm{d}_{s, N}^{i+1}(\mu)$ at timestep $t^{i+1}$ and for $\mu\in\mathcal{P}$:
\begin{gather}
\label{z}
\bm{z}_{f, N}^{i+1}({\mu}) = \sum_{k=1}^{N_z}\underline{z}_k^{i+1}({\mu})\bm{\Phi}_{\bm{z}_f}^k,\\
\label{p}
p_{f, N}^{0, i+1}({\mu}) = \sum_{k=1}^{N_{p}}\underline{p}_k^{0, i+1}({\mu})\Phi_{p}^k,\\
\label{d}
\bm{d}_{s, N}^{i+1}({\mu}) = \sum_{k=1}^{N_d}\underline{d}_k^{i+1}({\mu})\bm{\Phi}_{\bm{d}_s}^k.
\end{gather}
The reduced problem then reads: for every $i=0, \dots, N_T$ and for ${\mu}\in\mathcal{P}$: 
\subsubsection*{\underline{Mesh displacement}:}
let $\bm{d}_{f, N}^{i+1}({\mu})$ be defined by the reduced solid displacement at the previous time--step:
\begin{equation}
\label{mesh displacement reduced}
\bm{d}_{f, N}^{i+1}({\mu}) = \sum_{k=1}^{N_d}\underline{d}_k^i({\mu})\bm{\Phi}_{\bm{d}_f}^k;
\end{equation}
\subsubsection*{\underline{Fluid explicit step (with change of variable)}:}
find $\bm{z}_{f, N}^{i+1}({\mu})\in V_N$ such that $\forall \bm{v}_N\in V_N$:
\begin{equation*}
\begin{split}
&\rho_f\int_{\Omega_f} J_{\mu}K_{\mu_g}\Bigl(\frac{\bm{z}_{f, N}^{i+1}({\mu})-\bm{u}_{f, N}^i({\mu})}{\Delta T}\Bigr)\cdot\bm{v}_N\,dx + \mu_f\int_{\Omega_f}J_{\mu}K_{\mu_g}\varepsilon(\bm{z}_{f, N}^{i+1}({\mu}))\bm{F}^{-T}_{\mu}\bm{G}^{-T}_{\mu_g}:\nabla\bm{v}_N\,dx \\
& + \rho_f\int_{\Omega_f}J_{\mu}K_{\mu_g}\nabla\bm{z}_{f, N}^{i+1}({\mu})\bm{G}^{-1}_{\mu_g}\bm{F}^{-1}_{\mu}\bm{z}_{f, N}^{i+1}(\mu_g)\cdot\bm{v}_N\,dx  +\\
& +\rho_f\int_{\Omega_f}J_{\mu}K_{\mu_g}\nabla D_t\bm{d}_{f, N}^{i+1}({\mu})\bm{G}^{-1}_{\mu_g}\bm{F}^{-1}_{\mu}\bm{z}_{f, N}^{i+1}(\mu)\cdot\bm{v}_N\,dx \\
&+ \int_{\Omega_f}J_{\mu}K_{\mu_g}\bm{F}^{-T}_{\mu}\bm{G}^{-T}_{\mu_g}\nabla p_{f, N}^i({\mu})\cdot\bm{v}_h\,dx=-\rho_f\int_{\Omega_f} J_{\mu}K_{\mu_g}\Bigl(\frac{D_t\bm{d}_{f, N}^{i+1}(\mu)}{\Delta T}\Bigr)\cdot\bm{v}_N\,dx\\ 
&- \mu_f\int_{\Omega_f}J_{\mu}K_{\mu_g}\varepsilon(D_t\bm{d}_{f, N}^{i+1}({\mu}))\bm{F}^{-T}_{\mu}\bm{G}^{-T}_{\mu_g}:\nabla\bm{v}_N\,dx
 \quad\text{in $\Omega_f$.}
\end{split}
\end{equation*}
We then restore the reduced fluid velocity: $\bm{u}_{f, N}^{i+1}({\mu}) = \bm{z}_{f, N}^{i+1}({\mu}) + D_t\bm{d}_{f, N}^{i+1}({\mu})$.
\subsubsection*{\underline{Implicit step}:} for any $j=0, \dots$ until convergence:
\begin{enumerate}
\item \textbf{fluid projection substep:} find $p_{f, N}^{0, i+1, j+1}({\mu})\in Q_N^0$ such that $\forall q_N\in Q_N^0$:
\begin{equation*}
\begin{split}
& \alpha_{ROB}\int_{\Gamma_{FSI}}p_{f, N}^{0, i+1, j+1}(\mu)q_N\,ds +\\
& + \int_{\Omega_f}J_{\mu}K_{\mu_g}\bm{F}^{-T}_{\mu}\bm{G}^{-T}_{\mu_g}\nabla p_{f, N}^{0, i+1, j+1}(\mu)\cdot \bm{F}^{-T}_{\mu}\bm{G}^{-T}_{\mu_g}\nabla q_N\,dx=\\
&-\frac{\rho_f}{\Delta T}\int_{\Omega_f}\text{div}(J_{\mu}K_{\mu_g}\bm{G}^{-1}_{\mu_g}\bm{F}^{-1}_{\mu}\bm{u}_{f, N}^{i+1}(\mu))q_N\,dx + \alpha_{ROB}\int_{\Gamma_{FSI}}p_{f, N}^{0, i+1, j}(\mu)q_N\,ds -  \\
& -\alpha_{ROB}\int_{\Gamma_{FSI}}\ell_N^{i+1}q_N\,ds - \rho_f\int_{\Gamma_{FSI}}D_{tt}\bm{d}_{s, N}^{i+1, j}(\mu)\cdot J_{\mu}K_{\mu_g}\bm{F}^{-T}_{\mu}\bm{G}^{-T}_{\mu_g}\bm{n}_fq_N\,ds - \\
& \int_{\Omega_f}J_{\mu}K_{\mu_g}\bm{F}^{-T}_{\mu}\bm{G}^{-T}_{\mu_g}\nabla \ell_N^{i+1}\cdot \bm{F}^{-T}_{\mu}\bm{G}^{-T}_{\mu_g}\nabla q_N\,dx
\end{split}
\end{equation*}
we then recover the reduced fluid pressure $p_{f, N}^{i+1, j+1} = p_{f, N}^{0, i+1, j+1} + \ell_N^{i+1}$.
\item \textbf{structure projection substep:} find $\bm{d}_{s, N}^{i+1, j+1}({\mu})\in E^s_N$ such that $\forall \bm{e}_s\in E^s_N$:
\begin{equation*}
\begin{split}
&\rho_s\int_{\Omega_s}K_{\mu_g}D_{tt}\bm{d}_{s, N}^{i+1, j+1}(\mu)\cdot\bm{e}_N\,dx + \int_{\Omega_s}K_{\mu_g}\bm{P}(\bm{d}_{s, N}^{i+1, j+1}(\mu), \mu_s)\bm{G}^{-T}_{\mu_g}:\nabla\bm{e}_N\,dx = \\
&= -\int_{\Omega_s}J_{\mu}K_{\mu_g}\sigma_f(\bm{u}_{f, N}^{i+1}({\mu}), p_{f, N}^{i+1, j+1}({\mu}))\bm{F}^{-T}_{\mu}\bm{G}^{-T}_{\mu_g}\bm{n}_f\cdot\bm{e}_N\,dx.
\end{split}
\end{equation*}
\end{enumerate}

\subsection{Numerical results: geometrical parametrization only}\label{numerical results parametrized}
We now present some numerical results concerning the parametrized version of the two dimensional FSI test case presented in Section \ref{problem formulation}. The original domain is shown in Figure \ref{domain parametrization}, together with the reference configuration, and the parametrized reference configuration. The fluid domain is represented in blue, while the solid (the leaflets) is depicted in red. The geometrical constants defining the physical domain are reported in Table \ref{table geometrical leaflets ALE}. Only one geometrical parameter is considered here: the length $\mu_g$ of the leaflets, where we have chosen {$\mu_g\in\mathcal{P}=[0.8, 1.0]$}. An affine mapping $T_{\mu_g}$ is chosen to deform the reference domain $\hat{\Omega}$, obtained for $\mu_g=1.0$ cm, to the parametrized configuration $\tilde{\Omega}(\mu_g)$: {such a map is computed analytically}. Top and bottom walls of the blue domain are rigid, thus both the displacement $\bm{d}_f$ and the fluid velocity $\bm{u}_f$ are set to zero. Homogeneous Neumann condition is imposed on $\bm{u}_f$ on the outlet; a pressure profile $p_{in}(t)$ is described at the inlet, where:
\begin{equation*}
p_{in}(t) = 
\begin{cases}
5- 5\text{cos}\Bigl(\frac{2\pi t}{T_{in}}\Bigr) \qquad \text{for }t\leq 0.025s\\
5 \qquad \text{for }t>0.025s,
\end{cases}
\end{equation*}
and $T_{in}=0.1$ s. Also in this case we set a tolerance of $\varepsilon = 10^{-6}$ as a stopping criterion for the subiterations between the pressure problem and the solid problem.

\begin{table}
\centering
\caption{Physical and geometrical constants and parameters, for the geometrically parametrized leaflets test case.}
\label{table geometrical leaflets ALE}
\begin{tabular}{|l c|}
\hline
\textbf{Physical constants} & \textbf{Value}\\
\hline
$\rho_s$ & $1.1$ g/cm$^3$\\
$\mu_s$ & $10^5$\\
$\lambda_s$ & $8\times 10^5$\\
$\rho_f$ & $1$ g/cm$^3$\\
$\mu_f$ & $0.035$ Poise\\
\hline
\textbf{Geometrical constants} & \textbf{Value}\\
\hline
channel length & $10$ cm\\
channel height & $2.5$ cm\\
leaflets' thickness & $0.2$ cm\\
\textbf{Geometrical parameter} & \textbf{Value}\\
$\mu_g$ & $[0.8, 1.0]$\\
FE displacement order & $1$\\
FE velocity order & $2$\\
FE pressure order & $1$\\
\hline
\end{tabular}
\end{table}

For the simulation, {we use the same mesh used for the previous test case}, we set $\Delta t=10^{-4}$, for a maximum number of time--steps $N_T=500$, thus $T = 0.05$s. Table \ref{table geometrical leaflets ALE} summarizes the details of the offline stage and of the FE discretization. {The number of parameter samples used during the offline phase to train the algorithm varies between $N_g=8$ (for a total of $4000$ snapshots generated) to $N_g = 16$ (for a total of $8000$ snapshots generated).}

\begin{figure}[h]
\centering
\begin{minipage}{0.4\textwidth}
\centering
\includegraphics[width=0.9\textwidth]{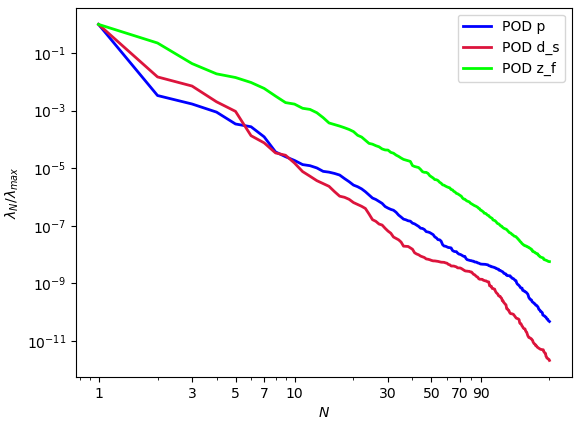}
\subcaption{POD eigenvalues}\label{1}
\end{minipage}%
\begin{minipage}{0.4\textwidth}
\centering
\includegraphics[width=0.9\textwidth]{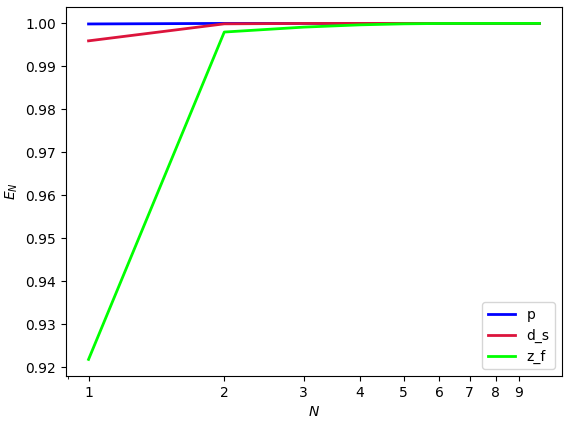}
\subcaption{POD retained energy.}\label{2}
\end{minipage}
\caption{POD results for the test case with a domain with geometrical parametrization.}
\end{figure}

\begin{figure}
\centering
\includegraphics[scale=0.3]{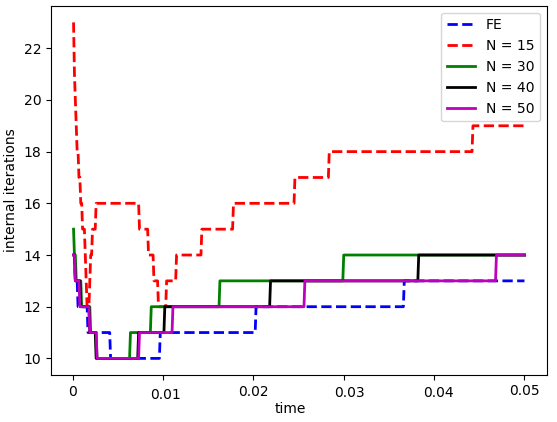}
\caption{Behavior in time of the number of iterations of the implicit step at the ROM level, compared to the FOM, with basis refinement. Number of parameter samples: $N_g=16$.}
\label{iterations in time geometrical}
\end{figure}

\begin{figure}
\centering
\includegraphics[scale=0.25]{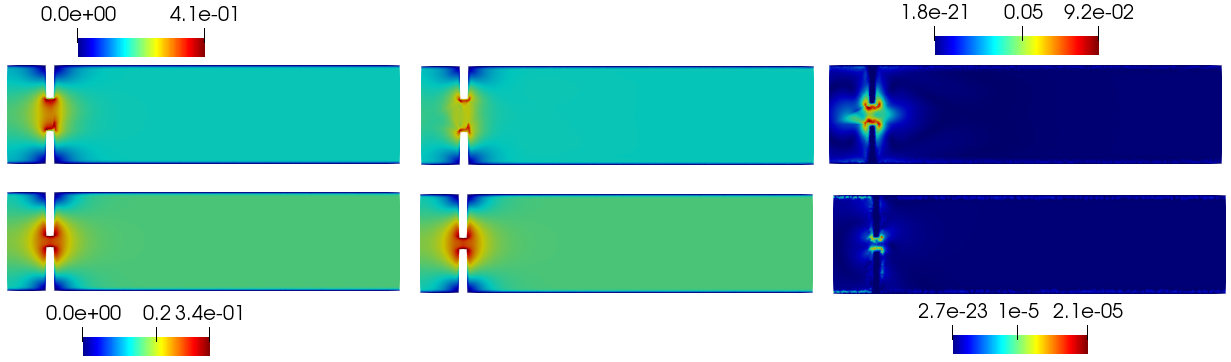}
\caption{{Representative solutions for the fluid velocity $\bm{u}_f$ at the final time $t=0.05$ s. Left column: Finite Element solutions for $\mu_g = 0.84$ (top) and $\mu_g = 1.0$ (bottom). Central column: reduced order solutions for $\mu_g = 0.84$ (top) and $\mu_g = 1.0$ (bottom). The reduced solutions were obtained with the reduced order model proposed, with $N=30$ basis for all the components. Right column: spatial distribution of the approximation error $\lvert \bm{u}_{f, h}-\bm{u}_{f, N}\rvert$ for $\mu_g = 0.84$ (top) and $\mu_g = 1.0$ (bottom).}}
\label{reduced velocity parametrized}
\end{figure}

\begin{figure}
\centering
\includegraphics[scale=0.4]{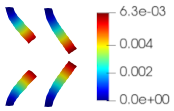}
\caption{Representative solutions for the displacement $\bm{d}_s$ at $t=0.05$ s, obtained with the reduced order model proposed ($N=30$ basis for all the components), for different values of the leaflet length $\mu_g$: $\mu_g = 0.84$ (left), and $\mu_g=1.0$ (right).}
\label{reduced displacement parametrized}
\end{figure}


\begin{table}
\centering
\caption{Average relative error of approximation for $\bm{u}_f$, $p_f$, $\bm{d}_s$, with sampling refinement. Leaflets' length $\mu_g = 0.84$ cm.}
\label{average error Ng geometric}
\begin{tabular}{|l c| c| c|}
\hline
$N_g$ & $\bm{u}_f$ & $p_f$ & $\bm{d}_s$\\
\hline
$9$ & $0.3225$ & $0.0208$ & $0.1662$\\
$13$ & $0.3205$ & $0.0210$ & $0.1636$\\
$16$ & $0.3245$ & $0.0204$ & $0.1679$\\
\hline
\end{tabular}
\end{table}

\begin{table}
\centering
\caption{Average relative error of approximation for $\bm{u}_f$, $p_f$, $\bm{d}_s$, with basis refinement. Leaflets' length $\mu_g = 0.84$ cm. Number of samples used: $N_g = 16$.}
\label{average error N geometric}
\begin{tabular}{|l c| c| c|}
\hline
$N$ & $\bm{u}_f$ & $p_f$ & $\bm{d}_s$\\
\hline
$15$ & $0.385$ & $0.0142$ & $0.1569$\\
$30$ & $0.3245$ & $0.020$ & $0.1679$\\
$40$ & $0.3138$ & $0.019$ & $0.1680$\\
$50$ & $0.2835$ & $0.0247$ & $0.1311$\\
\hline
\end{tabular}
\end{table}

\begin{figure}
\centering
\includegraphics[scale=0.3]{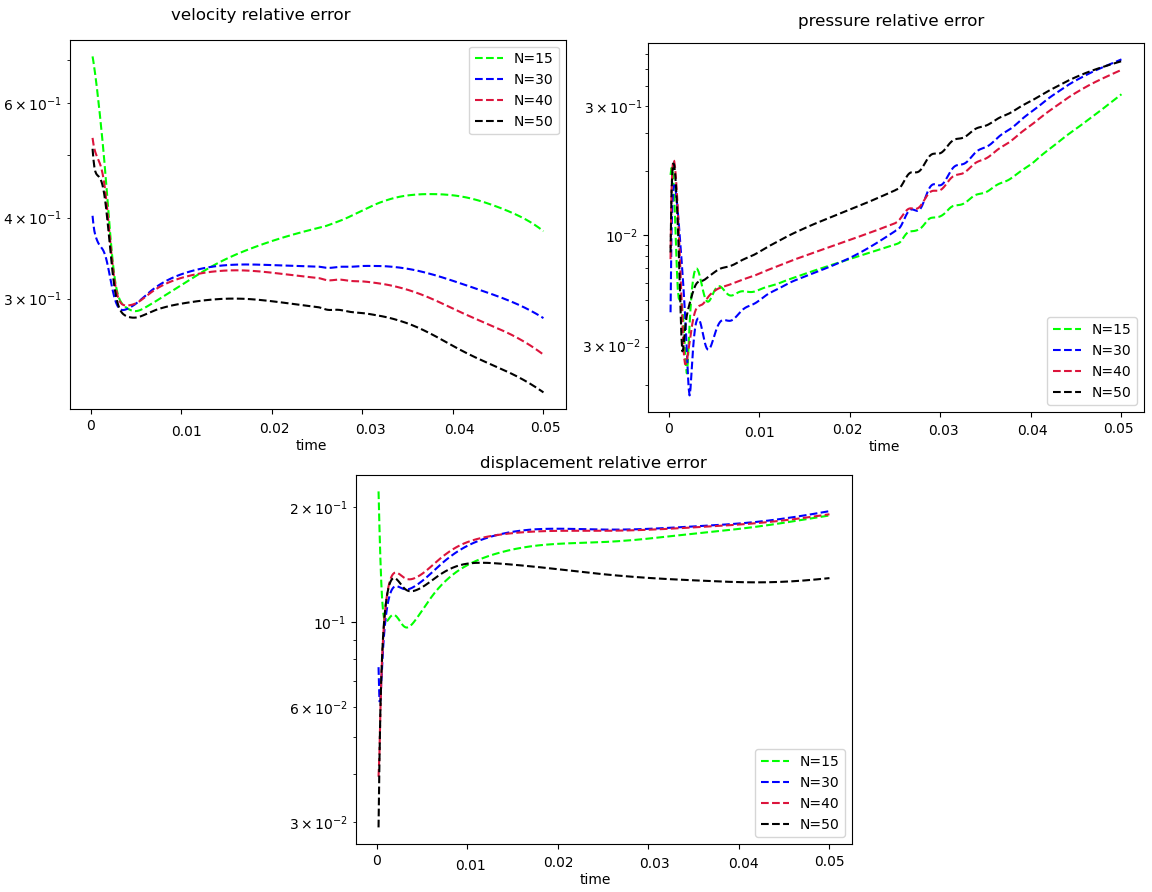}
\caption{{Relative error of approximation as a function of time. Top left: fluid velocity. Top right: fluid pressure. Bottom: solid displacement. The leaflets' length is $\mu_g=0.84$ cm. Number of sample parameters used: $N_g=16$.}}
\label{relative error in time}
\end{figure}

Figure \ref{1} shows the rate of decay of the eigenvalues returned by the final run of the POD on $\bm{z}_f$, $p_f$ and $\bm{d}_s$, respectively. As we can see, now the eigenvalues display and overall slower decay, showing that the complexity of the problem is higher with respect to the non parametrized case: this slower decay will then reflect into a higher number of modes that need to be used in the online phase. Figure \ref{2} shows the energy retained by the modes returned by the final run of the POD: the first mode of the pressure is indeed the most energetic one, while the first mode of the velocity is the least energetic one. All these results were obtained by using $N_g=16$ sampling parameters.
In Figure \ref{iterations in time geometrical} we depict the behavior in time of the number of iterations required in the implicit step, at the ROM level, to reach a tolerance of $\varepsilon=10^{-6}$, according to the number of modes used: for this test, the modes were generated using $N_g=16$ sampling parameters, for a total of $8000$ snapshots generated, and the leaflet length is $\mu_g=0.84$.  It is very interesting to see that, except for $N=15$, refining the number of basis functions used does not have a huge influence on the maximum number (maximum over the time interval under consideration) of internal sub--iterations, which is also quite close to the maximum number of subiterations required by the FOM.
{In Figure \ref{reduced velocity parametrized} we present a graphical comparison between two sets of reproductive solution for the fluid velocity, for two different values of the geometrical parameter $\mu_g=0.84$ and $\mu_g=1.0$, at the final timestep of the simulation $t=0.05$. The influence of the geometrical parameter is evident: the online solution represents very well the full order model solution, when $\mu_g=1.0$; on the other hand, when $\mu_g=0.84$, we can see that, with the number of modes used, we are not perfectly reproducing the correct magnitude of the jet between the two leaflets (the red area). This result highlights very well the deep influence that the geometrical parametrization has; indeed, the length $\mu_g=1.0$ corresponds to the reference leaflets length that we chose: in this case the geometrical deformation map $T_{\mu_g}$ defined in Section \ref{ALE geometrical formalism} is the identity. We would like to stress also the following: the test case with $\mu_g=0.84$ represents a prediction test case; indeed this value not only corresponds to a significant geometrical deformation, but it also corresponds to a value that has not been selected to train our algorithm at the offline level. For this reason, this test case is a stress test for our algorithm. In the framework of predictive problems, we believe that the result presented in Figure \ref{reduced velocity parametrized} can be improved, by refining the geometrical parameter sampling set used to train our method, thus by generating, during the offline phase, more snapshots: these simulations are very expensive, and for this reason we did not proceed any further with the sampling set refinement.} \\
Figure \ref{reduced displacement parametrized} shows the reduced displacement, for the same values of the geometrical parameter: $\mu_g=0.84$ and $\mu_g=1.0$; again, the influence of $\mu_g$ is clear: the longer the leaflets, the bigger their deformation is going to be, under the same physical parameters.
Table \ref{average error Ng geometric} represent the average approximation error for the fluid velocity, the fluid pressure and solid displacement, with refinement of the training sample $\mathcal{P}_{train}$: all the reduced solutions have been obtained using $N=30$ reduced basis functions for all the components $\bm{u}_f$, $p_f$ and $\bm{d}_s$; it is interesting to see that, for the highest number of training samples, namely $N_g=16$ (which corresponds to a total of $8000$ snapshots generated) we observe a slight increase in the average approximation error of the velocity and the solid displacement. We think this could be due to the fact that the online model could benefit from more reduced basis for $\bm{u}_f$ and $\bm{d}_s$, and our hypothesis seems to be confirmed by the results of Table \ref{average error N geometric}, where we show the average approximation error for $N_g=16$, with basis refinement: as we can see, the model benefits from the increment in the number of modes used in the online simulations.
Also here the average error is intended as average over time, and it is computed as the $H^1$ relative error for the velocity, the $L^2$ error for the pressure and the $L^2$ error for the displacement. {Finally, in Figure \ref{relative error in time} we present the behavior of the relative approximation error in time: we consider the reduced solution to be obtained with $N=15, 30, 40, 50$ modes for each component (see color legend in the figure). As we observe, the behavior in time of the approximation error seems to confirm the results reported in Table \ref{average error N geometric}: our model benefits from an increase in the number of modes used. It is interesting to observe the behavior of the pressure relative error: the error accumulates over time, thus steadily increases. This represents a starting point for future studies and future work development; indeed, this steady growth in time of the approximation error is related only to the pressure component of the FSI solution. We ask ourselves if this is somehow related to the fact that we are not using a supremizer enrichment of the velocity space at the online level; it would be therefore interesting to study if this has some effect in the online approximation of the pressure, especially for long time simulations.}
\subsection{Numerical results: physical and geometrical parametrization}\label{shape and physical parametrization: numerical results}
We present some numerical results for the test case with a geometrical and physical parametrization: now $\mu=(\mu_g, \mu_p)$, where the physical parameter $\mu_p$ represents the shear modulus of the leaflets, and thus $\mu_p=\mu_s$, the second Lam\'e constant. The original configuration, the intermediate configuration and the reference configuration are the ones depicted in Figure \ref{domain parametrization}. The geometrical constants of the problem are the same ones of the previous test case, and are reported in Table \ref{table geometrical leaflets ALE}. The physical parameter $\mu_s$ varies within the range $[10^5, 8\times 10^5]$. We use the same boundary conditions, the same inlet pressure profile, time step and tolerance used in the previous test case.

\begin{figure}[h]
\centering
\begin{minipage}[b]{0.45\textwidth}
 \centering
  \includegraphics[width=0.25\textwidth]{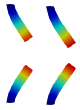}
   \subcaption{Same shear modulus $\mu_s=10^5$, but increased length $\mu_g = 0.82$ (left) and $\mu_g=1$ (right)}\label{length_long}
\end{minipage}%
\hfill
\vline
\begin{minipage}[b]{0.45\textwidth}
  \centering
  \includegraphics[width=0.5\textwidth]{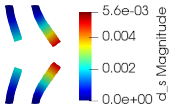}
 \subcaption{Same leaflet length $\mu_g=0.82$cm, but increased shear modulus $\mu_s=10^5$ (left) and $\mu_s=8\times 10^5$ (right)}\label{length_short}
\end{minipage}%
\hfill
\caption{Reduced order solid displacement $\bm{d}_s^N(\mu)$. Comparison of different behaviors of the material, for different values of the geometrical and physical parameters. From left to right: same leaflets length (length of $0.8$cm) and increased shear modulus ($\mu_s=100000$, $800000$); same leaflet length (length of $1$cm) and increased shear modulus ($\mu_s=100000$, $800000$); increased leaflet length ($\mu_g=0.8$, $1.0$ cm), and same shear modulus $\mu_s=100000$.}\label{comparison}
\end{figure}

\begin{figure}
\centering
\includegraphics[scale=0.3]{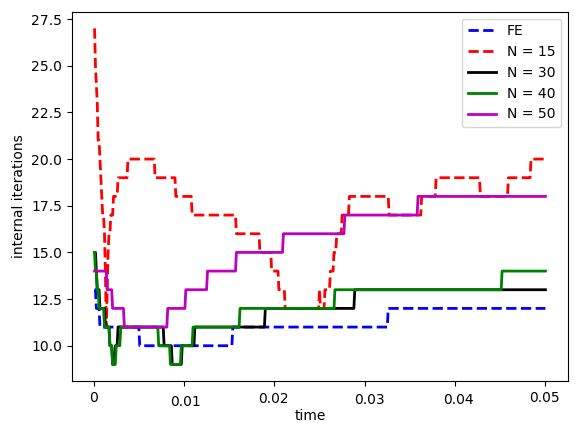}
\caption{Number of subiterations for the implicit online step, as a function of time, according to the number of reduced basis used online.}
\label{internal iterations in time physical parametrization}
\end{figure}

\begin{figure}
\centering
\includegraphics[scale=0.3]{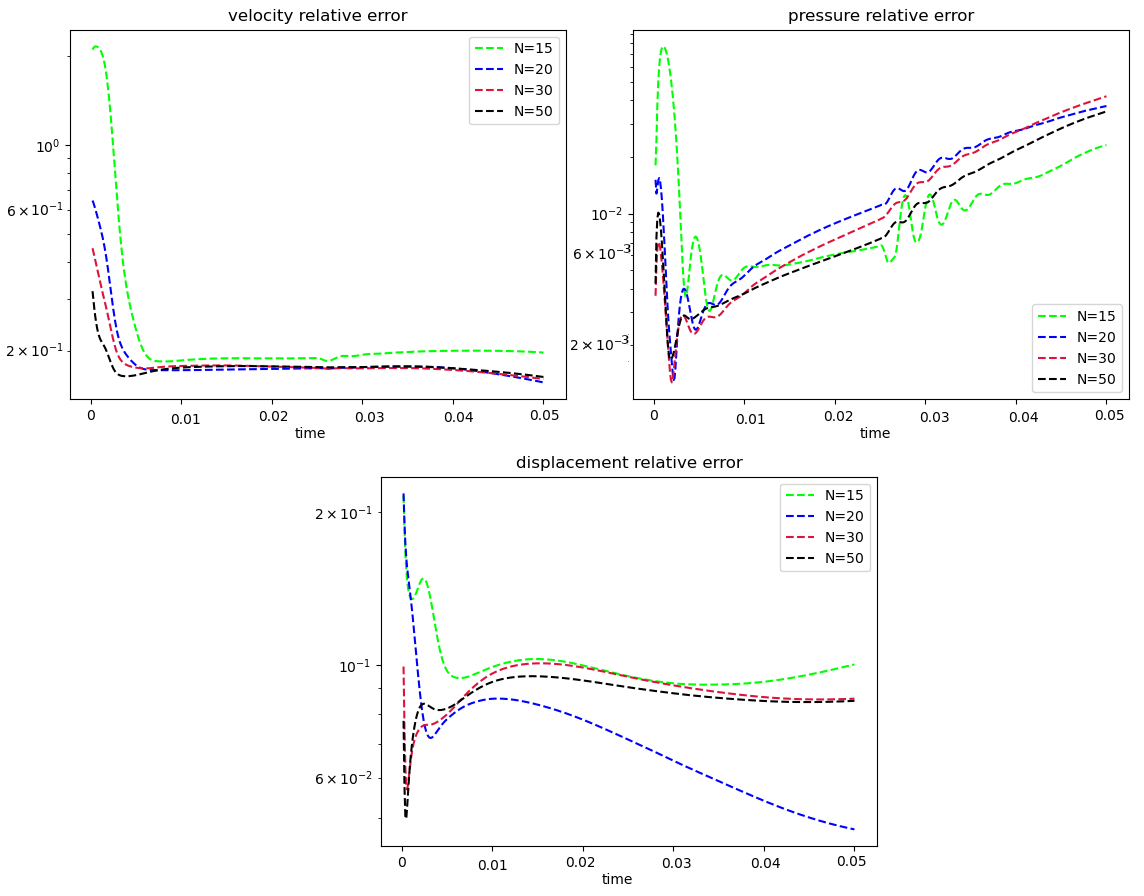}
\caption{{Behavior in time of the relative approximation error, for different number of modes $N$ used in the online model. Top left: velocity error. Top right: pressure error. Bottom center: displacement. Leaflets' length $\mu_g=0.9$ cm, leaflets' shear modulus $\mu_s= 10^5$. Number of sampling parameters used: $N_g = 8$, $N_s = 10$.}}
\label{errors in time physical parametrization}
\end{figure}

\begin{table}
\centering
\caption{Average relative error of approximation for $\bm{u}_f$, $p_f$, $\bm{d}_s$, with basis refinement. Number of parameters samplings: $N_g = 8$ and $N_s=10$. Leaflets' length $\mu_g = 0.9$ cm and shear modulus $\mu_s = 10^5$.}
\label{average error N}
\begin{tabular}{|l c| c| c|}
\hline
$N$ & $\bm{u}_f$ & $p_f$ & $\bm{d}_s$\\
\hline
$15$ & $0.28241$ & $0.01206$ & $0.09978$\\
$20$ & $0.18863$ & $0.01562$ & $0.07082$\\
$30$ & $0.18152$ & $0.01489$ & $0.09008$\\
$50$ & $0.17635$ & $0.01225$ & $0.08762$\\
\hline
\end{tabular}
\end{table}

\begin{table}
\centering
\caption{Average relative error of approximation for $\bm{u}_f$, $p_f$ and $\bm{d}_s$, with refinement of the geometrical parameter sampling. Number of physical parameter samplings: $N_s=10$. Leaflets' length $\mu_g = 0.9$ cm and shear modulus $\mu_s = 10^5$.}
\label{average error Ng}
\begin{tabular}{|l c| c| c|}
\hline
$N_g$ & $\bm{u}_f$ & $p_f$ & $\bm{d}_s$\\
\hline
$5$ & $0.18206$ & $0.01419$ & $0.09073$\\
$6$ & $0.18108$ & $0.01460$ & $0.09054$\\
$8$ & $0.18152$ & $0.01489$ & $0.09008$\\
\hline
\end{tabular}
\end{table}

\begin{table}
\centering
\caption{Average relative error of approximation for $\bm{u}_f$, $p_f$ and $\bm{d}_s$, with refinement of the physical parameter sampling. Number of geometrical parameter samplings: $N_g=8$. Leaflets' length $\mu_g = 0.9$ cm and shear modulus $\mu_s = 10^5$.}
\label{average error Ns}
\begin{tabular}{|l c| c| c|}
\hline
$N_s$ & $\bm{u}_f$ & $p_f$ & $\bm{d}_s$\\
\hline
$6$ & $0.19247$ & $0.01151$ & $0.1200$\\
$7$ & $0.19121$ & $0.01174$ & $0.1185$\\
$9$ & $0.18415$ & $0.01147$ & $0.1108$\\
$10$ & $0.18152$ & $0.01489$ & $0.0900$\\
\hline
\end{tabular}
\end{table}

\begin{table}
\centering
\caption{CPU time comparison and average approximation error, for three different tolerances $\varepsilon$ of the implicit step at the ROM level. {The FOM CPU time is $24763.47$ seconds for $\varepsilon = 10^{-6}$.}}
\label{CPU times}
\begin{tabular}{|l c| l| c| c| c|}
\hline
$\varepsilon$ & ROM CPU& avg.error $\bm{u}_f$ & avg.error $p_f$ & avg.error $\bm{d}_s$\\
\hline
$10^{-3}$ & $9979.5 $& $0.16309$ & $0.03300$ & $0.11955$\\
$10^{-6}$ & $13696$ & $0.16302$ & $0.03296$ & $0.11933$\\
$10^{-8}$ & $16088.5$ & $0.16302$ & $0.03296$& $0.11933$\\
\hline
\end{tabular}
\end{table}

Figure \ref{comparison} shows two different examples of the behavior of the leaflets, according to the change of the physical and/or of the geometrical parameters: all the pictures represent the online displacement of the solid at the final timestep of the simulation, namely for $t=0.05$s. The results have been obtained using $N=30$ reduced basis for all the components; as we can see from Figure \ref{length_short} and \ref{length_long}, an increase in the shear modulus leads to a material that is much more hard to deform. On the contrary, for a fixed value of the properties of the material under consideration, and increase of the length of the leaflets leads to an increase in the displacement. 
{Figure \ref{internal iterations in time physical parametrization} shows the behavior in time of the number of iterations of the implicit step, for different number of modes $N$ used in the online phase, compared against the FOM. Figure \ref{errors in time physical parametrization} represent the relative error approximation for $\bm{u_f}$ (top left), $p_f$ (top right) and $\bm{d}_s$ (bottom center) as a function of time. From Figure \ref{errors in time physical parametrization} we can see that the relative error for the fluid velocity stabilizes after some iterations, reaching a magnitude of $2\times10^{-1}$ when using $N= 20$, $30$ and $40$ modes. Also the relative error for $\bm{d}_s$ shows a plateau around $9\times10^{-2}$, except for when $N=20$, in which case the error decreases in time: we read this result as the consequence of using too many modes in the online model, when $N=30, 50$. For the fluid pressure we observe the same accumulation phenomenon that we observe in the previous parametrized test case. To conclude, even though these relative errors seem high, we are again testing our algorithm with a prediction problem, since the value $(\mu_g, \mu_s)=(0.9, 10^5)$ has not been used to generate modes during the training phase of the algorithm. In addition to this, we would like to remark that, by increasing the number of sampling parameters used at the FOM level, we should be able to drive the error down, and this seems to be confirmed by the results in Table \ref{average error Ns}: we did not proceed with a sampling refinement because of time constraints, since the generation of $40,000$ snapshots is very demanding.
}
{In Tables \ref{average error N}, \ref{average error Ng} and \ref{average error Ns} we report the average approximation error for the fluid velocity $\bm{u}_f$, the fluid pressure $p_f$ and the solid displacement $\bm{d}_s$, when $\mu_g = 0.9$ cm and $\mu_s=10^5$: the average has been taken over time, and we used the $L^2$ norm for the fluid pressure and the solid displacement, and the $H^1$ norm for the fluid velocity. In Table \ref{average error N} we computed the relative approximation error, by using $N_g=8$ and $N_s=10$ training samples, and refining the number of basis functions used online, from $N=15$ to $N=50$. In this case, we used the same number of basis functions for all the components of the FSI solution: again, we remark that other tests are possible, for the interested reader, for example keeping the number of modes fixed for the fluid velocity and the solid displacement, and varying the number of modes used for the fluid pressure. Also here, like for the previous test case, the online model seems to benefit from a higher number of modes used. In Table \ref{average error Ng}, we used $N=30$ modes for the online solution, $N_s=10$ training samples for the physical parameter, and we refined the geometrical training set, from $N_g=5$ to $N_g=8$. On the contrary, in Table \ref{average error Ns}, we kept $N_g=8$ training samples for the geometrical parameter, and refined the physical training set from $N_s=6$ to $N_s=10$. We remark that analyzing the training samples with $N_g=8$ and $N_s=10$ corresponds to generating $40.000$ snapshots. Indeed, in the framework of the POD, we have to compute all the snapshots for each value of the training parameter, in an unsteady framework: this procedure is extremely expensive, and it required a total of $5$ days, by running the simulations on two different computers: a computer Intel(R) Core(TM) i5-4670S CPU with 3.10 GHz and 7.7G of RAM, and a supercomputer Intel(R) Xeon(R) CPU E5-2687W v4 with 3.00 GHz and 540G of RAM.
For this reason, we did not proceed further with the refinement of the parameter sampling. As the results show, however, by refining the parameter space, thus by using more training samples, we are able to improve the average approximation error. We remark also that, in case of Table \ref{average error Ng} and Table \ref{average error Ns}, the online computations have been made using $N=30$ modes for each component of the solution. We think that in this case the average approximation errors are higher (compared to the ones obtained for the non parametrized test case) expecially due to the presence of the geometrical parametrization. 
Finally, in Table \ref{CPU times} we show some CPU times: we choose $\mu_g=0.9$ and $\mu_s=8\times10^5$, and fix $N=30$ for all the components (these modes were obtained fixing $N_g=8$ and $N_s=10$). The simulation is run on a computer Intel(R) Core(TM) i5-4670S CPU with 3.10 GHz and 7.7G of RAM; the CPU execution time is reported in Table \ref{CPU times} (times measured in seconds): the FOM requires $24763.47$ seconds for the same geometrical and physical parameters (the computation comprises also the computation of the change of variable $\bm{z}_f$ and the homogenization of the fluid pressure through the lifting function).
As we can see, and as expected, by strengthening the tolerance, the computational time grows, as more iterations are needed at the implicit step to reach the convergence. What is interesting to notice, however, is that,  by strengthening the tolerance we do not see an important decrease in the average approximation error: this result seems to therefore suggest that online computations could be carried out also using a coarser tolerance (for example $10^{-4}$ instead of $10^{-6}$).
}

\section{Conclusions}\label{conclusions}
In this manuscript we presented a Reduced Order Model algorithm designed to address FSI problems, in the unsteady case, and possibly in the presence of a parameter dependence. The ROM is based on a partitioned procedure: the main advantage of this is given by the fact that, by solving separately the fluid and the solid problem, we are not only able to lower the dimension of the systems to be solved in the online phase, but we also have a better control on the number of variables that are needed. 
The reduced basis functions are generated through a Proper Orthogonal Decomposition on the set of snapshots, with the introduction of a change of variable in the fluid problem formulation.
The procedure that we have proposed aims at extending the work presented in \cite{Ballarin2017, monica} to the case of the coupling between an incompressible fluid and a thick, two dimensional structure, also in the presence of geometrical parametrization. The results that we have obtained confirm the following aspects:
\begin{itemize}
\item introducing a change of variable in fluid explicit step allows us to avoid the introduction of a further unknown in the system, namely a Lagrange multiplier, in order to impose non homogeneous boundary conditions at the fluid--structure interface;
\item the choice of not performing a POD on the snapshots $\bm{d}_f$, but rather  performing an harmonic extension of some modes, allows us to build the online $\bm{d}_{N, f}$ in a cheap way. 
Moreover, the coupling condition that imposes the continuity of the displacements at the fluid--structure interface is automatically satisfied, thanks to the way we have defined the reduced basis for $\bm{d}_f$.
\end{itemize} 
In addition to the list of remarks presented, another very important detail of the procedure presented in this manuscript is the following: we did not rely a supremizer enrichment of the fluid velocity space, as it is the usual case for reduction methods, in order to obtain a stable approximation of the fluid pressure in the online phase. Our choice is motivated by the fact that, even at the Finite Element level, the Chorin--Temam projection scheme with the pressure Poisson formulation can be applied succesfully also to velocity--pressure FE spaces that do not satisfy the inf--sup condition, see \cite{GuermondQuartapelle2}. This allows us to limit the dimension of the system to be solved online, and it is a big motivation for the choice of a Chorin--Temam projection scheme within our partitioned approach.\\
While testing this algorithm, we have seen that a partitioned procedure is demanding from the computational time point of view: this drawback is represented by the fact that, in the imposition of the coupling conditions through a Robin boundary condition, the constant $\alpha_{ROB}$ that makes the procedure more stable depends on the time--step used. {If we choose a time--step that is too big, then our Robin coefficient $\alpha_{ROB}$ becomes very small, and we recover the original Dirichlet--Neumann coupling, which is known to have stability problems, i.e. the implicit step may not converge. Always in the direction of the computational effort of the offline phase, we remark that in this manuscript we considered a physical parameter for the solid, but no parametrization of the fluid has been taken into account: considering a fluid parameter, instead of a solid one, does not change the design of the algorithm. However having a test case which comprises both a parameter for the solid and a parameter for the fluid does highlight the boundaries of the POD in this framework; indeed, generating the snapshots for the chosen sampling set would be extremely time demanding, even though it would be carried out at an offline level: for this reason, future perspectives also include the design of an error estimator and the integration of a Greedy algorithm based on such estimator.}
{It is also important to mention the fact an efficient online--offline decoupling is very important in terms of model order reduction efficiency: for the fluid part, this can be recovered thanks to the Empirical Interpolation Method (EIM), see for example \cite{BARRAULT2004667, MadayMulaPateraYano, BallarinRozza}, whereas for the solid mechanics part, one may think about using hyper--reduction procedures that preserve the Hamiltonian structure, such as, for example, the ECSW \cite{farhatECWS}. These techniques have not been used in this work, as we wanted to focus on the development and test of a reduced order segregated procedure for FSI problems which involve the coupling of an incompressible fluid with an elastic structure; this can be seen as a natural future perspective for this work.\\} {As it was mentioned in the Introduction, there is currently a fair amount of interest in approaches that are able to couple high fidelity models with reduced order models: one may think about using a FOM for the structure, and a ROM for the fluid. The authors believe that there is a possible extension (with some required modifications) of the algorithm proposed, to this kind of situations: again, this represents another interesting future line of research.
Another remark that we would like to make concerns the timestep used in our numerical simulations: in the results that we have showed, a timestep of $\Delta T=10^{-4}$ has been used both for the fluid and for the solid problem. However, one may be interested in using different timesteps for the two physics, as the solid model is an hyperbolic equation, and the fluid is a parabolic one: we do believe that this is possible, with the procedure presented. Indeed, given the time interval $I_n:=(t_{n-1}, t_n]$, given a solid timestep $\Delta T_s$ and a fluid timestep $\Delta T_f$, assume that the resulting time discretization $I_n^s$ of the interval $I_n$ is finer than the time discretization $I_n^f$ of  $I_n$. Then, one just needs to be able to evaluate fluid quantities on the times $t_i^s\in I_n^s$ and solid quantities on the times $t_i^f\in I_n^f$: with the definition of suitable interpolation operators (interpolation in time), we should be able to implement a partitioned scheme with different timesteps for the two physics. This idea is presented for example in \cite{RichterSoszynska} for a monolithic scheme, but it represents an interesting starting point for a future work within partitioned schemes. In this case, maybe manually tuning the Robin parameter $\alpha_{ROB}$ to achieve optimal convergence can be a better idea, in order to drop the dependence on $\Delta T$: further research in this direction needs to be carried out.
}
\section*{Acknowledgements}
We acknowledge the support by European Union Funding for Research and Innovation -- Horizon 2020 Program -- in the framework of European Research Council Executive Agency: Consolidator Grant H2020 ERC CoG 2015 AROMA-CFD project 681447 ``Advanced Reduced Order Methods with Applications in Computational Fluid Dynamics'' (PI Prof. Gianluigi Rozza).
We also acknowledge the INDAM-GNCS project ``Tecniche Numeriche Avanzate per Applicazioni Industriali''. The author M. Nonino also acknowledges the support of the Austrian Science Fund (FWF) project F65 ``Taming complexity in Partial Differential Systems'' and the Austrian Science Fund (FWF) project P 33477.
Numerical simulations have been obtained with the extension \emph{multiphenics} of \emph{FEniCS} \cite{multiphenics, fenics} for the high fidelity part, and RBniCS \cite{rbnics} for the reduced order part. We acknowledge developers and contributors of each of the aforementioned libraries.

\bibliographystyle{abbrv}
\bibliography{main}

\end{document}